\documentstyle{amsppt}
\magnification 1200
\UseAMSsymbols
\hsize 6.0 true in
\vsize 8.5 true in
\parskip=\medskipamount
\NoBlackBoxes

\def\mathcal{\Cal}

\def\supp{\text{\rm supp\,}}

\def\ve{\varepsilon}
\def\vp{\varphi}

\def\snint{\raise2pt\hbox{$_{^\not}$}\kern-3.5 pt\int}
\def\nint{\not\!\!\int}
\def\dist{\text{\rm dist\,}}

\def\snint{\raise2pt\hbox{$_{^\not}$}\kern-3.5 pt\int}

\def\mes{\text {\rm mes\,}}
\TagsOnRight
\NoRunningHeads

\document
\parskip=\medskipamount
\topmatter
\title
On the Schr\"odinger maximal function in higher dimension
\endtitle
\author
J.~Bourgain\\
\\
(dedicated to B.~Kasin for his 60$^{\text{th}}$)
\endauthor
\address
Institute for Advanced Study, Princeton, NJ 08540
\endaddress
\email
bourgain\@ias.edu
\endemail
\abstract
New estimates on the maximal function associated to the linear Schr\"odinger equation are established.
\endabstract
\endtopmatter

\bigskip
\noindent
{\bf 1. Introduction}

Recall that the solution of the linear Schr\"odinger equation
$$
\left\{
\aligned
&iu_t-\Delta u =0\\
&u(x, 0) =f(x)
\endaligned
\right.
\tag 1.0
$$
with $(x, t)\in \Bbb R^n\times \Bbb R$ is given by
$$
e^{it\Delta} f(x) =(2\pi)^{-n} \int e^{i(x. \xi + t|\xi|^2)} \widehat f(\xi) d\xi.\tag 1.1
$$
Assuming $f\in H^s (\Bbb R^n)$ for suitable $s$, when does the almost convergence property
$$
\lim_{t\to 0} e^{it\Delta} f = f \ a.e. \tag 1.2
$$
hold?
This problem originates from Carleson's paper \cite {C} who proved convergence for $s\geq\frac 14$ when $n=1$.
Dahlberg and Kenig \cite{D-K} showed that this result is sharp.  
In dimensions $n\geq 2$, Sh\"olin \cite{S} and Vega {V} established independently convergence for $s>\frac 12$, while similar examples
as considered in \cite {D-K} show failure of convergence for $s<\frac 14$.

The problem for $n=2$ has been studied by various authors.
Proof of convergence for some $s<\frac 12$ appears first in the author's papers \cite {B1}, \cite{B2}.
Subsequently, the required threshold was lowered by Moyua, Vargas, Vega \cite {M-V-V}, Tao, Vargas \hbox{\cite{T-V1,2}} and S.~Lee \cite{L}.
The strongest result to date appears in \cite {L} and asserts a.e. convergence for $f\in H^s(\Bbb R^2), s>\frac 38$.
On the other hand, no improvements of Sj\"olin-Vega result for $n\geq 3$ seem to have been obtained so far.

We will prove the following

\proclaim
{Theorem 1}  
For all $n$, there is $\ve_n>0$ such that the a.e. convergence property holds for $f\in H^s(\Bbb R^n), s>\frac 12-\ve_n$.

In fact, one may take $\ve_n =\frac 1{4n}$.
\endproclaim

Note that for $n=2$, the exponent corresponds to the one obtained in \cite {L}.
The reason for formulating Theorem 1 this way is that the first (qualitative) statement allows for a less involved argument that will be
presented first.

Most of the progress on this problem relates to the advances around the restriction theory for the paraboloid.
In particular, the method used in \cite {L} depends on Tao's bilinear restriction estimate \cite {T}.
The proof of Theorem 1 will be based on the results and techniques from \cite {B-G} and hence ultimately on the multi-linear restriction theory
developed in \cite {B-C-T}.

As pointed out in \cite {L}, it should also be observed that all results obtained in the convergence problem for (1.0) hold equally well for
generalized Schr\"odinger equations
$$
iu_t -\Phi (D) u=0\tag 1.3
$$
with $\Phi(D)$ a Fourier multiplier operator satisfying $|D^\alpha \phi (\xi)|\leq C|\xi|^{2-\alpha}, \break 
|\nabla\phi(\xi)|\geq c|\xi|$.
This is also the case for Theorem 1.

Perhaps the most interesting point in this Note is a disproof of what one seemed to believe, namely that $f\in H^s(\Bbb R^n), s>\frac 14$, should be
the correct condition in arbitrary dimension $n$.

\proclaim
{Theorem 2}  For $n>4$, the a.e. convergence property for the Schr\"odinger group requires $f\in H^s(\Bbb R^n) $ with $s\geq \frac {n-2}{2n}$.
\endproclaim

Hence this exponent $\to \frac 12$ for $n\to\infty$.
The examples are of a quite different nature from previous constructions and rely on an arithmetical input, namely lattice prints on spheres.

Returning to the generalized setting, we exhibit in the last section of the paper some further examples setting further restrictions to what may be
proven in the generality of our approach to Theorem 1.

\bigskip

\noindent
{\bf 2. Braking the $\frac 12$-barrier}

Before getting to more quantitative analysis, we start by establishing the first part of Theorem 1.
The main input in that argument is the multi-linear bound from \cite {B-C-T} and it is technically rather easy.

A few preliminary reductions.
Assume $R>1$ large and
$$
\supp \widehat f \subset [|\xi|\sim R].\tag 2.0
$$
Our aim is to establish a maximal inequality of the form
$$
\big\Vert \sup_{0<t<1}| e^{it\Delta} f|\big\Vert_{L^2(B_1)}\leq CR^\alpha \Vert f\Vert_2\tag 2.1
$$
where $B_K=[x\in\Bbb R^n; |x|\leq K]$.
From the stationary phase, (2.1) may be deduced from the restricted inequality
$$
\big\Vert \sup_{0<t<\frac 1R} |e^{it\Delta} f|\big\Vert_{L^2(B_1)} \leq CR^{\alpha'} \Vert f\Vert_2\tag 2.2
$$
for any $\alpha' <\alpha$ (see \cite{L}, Lemma 2.3).
Hence we restrict $0<t<\frac 1R$ in (2.1).

Next, performing a rescaling, we need to show that
$$
\big\Vert \sup_{0<t<R}|e^{it\Delta} f|\big\Vert_{L^2(B_R)} \leq CR^\alpha \Vert f\Vert_2\tag 2.3
$$
assuming
$$
\supp\widehat f  \subset[|\xi|\sim 1].\tag 2.4
$$
This format corresponds to \cite{B-G}.
What follows is closely related to arguments in \S2, \S3 of that paper.

For notational convenience, set $x_{n+1} =t$ and
$$
\psi (x, \xi)= x_1 \xi_1 +\cdots+ x_n\xi_n+x_{n+1} |\xi|^2.\tag 2.5
$$
Hence
$$
(e^{it\Delta} f)(x_1, \ldots, x_n)=\int_\Omega\widehat f (\xi) e^{i\psi(x, \xi)}d\xi \equiv T\widehat f\tag 2.6
$$
with
$$
\Omega =[|\xi|\sim 1].
$$
Fix $1\ll K= K (R) \ll R$ to be specified later.
Partition $\Omega$ in boxes $\Omega_\alpha$ of size $\frac 1K$ centered at $\xi_\alpha$. Then
$$
Tf=\sum_\alpha e^{i\psi(\alpha, \xi_\alpha)} T_\alpha f\tag 2.7
$$
with
$$
T_\alpha f(x) =\int_{\Omega_\alpha} f(\xi) e^{i[\psi(x, \xi)-\psi(x, \xi_\alpha)]} d\xi.\tag 2.8
$$
Note that $|\nabla_x [\psi(x, \xi)-\psi(x, \xi_\alpha)] \lesssim \frac 1K$ and hence $T_\alpha f$ may be viewed as `approximately constant' 
on balls of size $O(K)$ in $B_R\subset\Bbb R^{n+1}$.
We refer the reader to \cite {B-G}, \S2, for the technical details needed to formalize this last statement.

Fix a ball $B_K\subset B_R$.
The main idea is to exploit a dichotomy similar to the one used in \S2, \S3 of \cite {B-G}.
More precisely, we distinguish the following two alternatives.

\noindent
{\bf Case I: `$(n-1)$-transversality'}

Viewing the $T_\alpha f$ as essentially constant on $\Omega_\alpha$, this means that there are indices 
$\alpha_1, \ldots, \alpha_{n+1}$ such that

\roster
\item "{(2.9)}" $|\nu(\xi_1')\wedge\ldots\wedge \nu(\xi'_{n+1})|>K^{-C}$ whenever $\xi_j'\in \Omega_{\alpha_j}$.

Here $\nu (\xi)// -2 \sum^n_{i=1} \xi_ie_i +e_{n+1}$ denotes the normal of the hypersurface $(\xi, |\xi|^2)$
\endroster

\noindent
and 

\medskip
\roster
\item "{(2.10)}" $|T_{\alpha_1} f|, \ldots, |T_{\alpha_{n+1}} f|\geq K^{-n} \max_\alpha |T_\alpha f|$.
\endroster
\medskip

Hence the following majoration holds on $B_K$
$$
|Tf|\lesssim K^n \max_\alpha |T_\alpha f|\lesssim K^{2n} \big[\prod^{n+1}_{j=1} |T_{\alpha_j} f|\big]^{\frac 1{n+1}}
\leq K^{2n}\sum_{\Sb \alpha_1, \ldots, \alpha_{n+1}\\ \text { transversal}\endSb } \Big[\prod^{n+1}_{j=1} |T_{\alpha_j} f|\Big]^{\frac 1{n+1}}.
$$
Recall the multi-linear restriction bound from \cite {B-C-T}, implying that for $\alpha_1, \ldots, \alpha_{n+1}$ `transversal' as defined in (2.9)
$$
\Big\Vert\Big[\prod^{n+1}_{j=1} |T_{\alpha_j} f|\Big]^{\frac 1{n+1}} \Big\Vert_{L^q(B_R)} \ll K^CR^\ve\Vert f\Vert_2\tag 2.12
$$
holds, with $q=\frac {2(n+1)}n$.

Hence, denoting $x=(x', x_{n+1})$ and using H\"older's inequality, the corresponding contribution to the maximal function may be bounded by
$$
\align
&\Vert \sup_{|x_{n+1}|<R} (2.11) \Vert_{L^2_{[|x'|<R]}}\lesssim\\
&\Vert (2.11) \Vert_{L^2_{[|x'|<R]}} L^q_{[|x_{n+1}|< R]}\lesssim\\
&R^{n(\frac 12-\frac 1q)} \Vert(2.11)\Vert_{L^q_{[|x|<R]}}\overset {(2.12)}\to \ll\\
&R^{n(\frac 12-\frac 1q)+\ve} K^C\Vert f\Vert_2 =R^{\frac n{2(n+1)}+\ve} K^C\Vert f\Vert_2.\tag 2.13
\endalign
$$

\medskip
\noindent
{\bf Case II: Failure of $(n+1)$-transversality}

In this situation, there is an $(n-1)$-dimensional affine hyperplane $\Cal L$ in $\Bbb R^n$ such that $\dist (\Omega_\alpha, \Cal L)\lesssim\frac
1K$ if $|T_\alpha f|\geq K^{-n}\max_\alpha |T_\alpha f|$ on $B_K$.

Denoting $\tilde {\Cal L}$ and $\frac 1K$-neighborhood of $\Cal L$, it follows that on $B_K$
$$
|Tf(x)|\leq \Big| \sum_{\Omega_\alpha\subset\tilde {\Cal L}} e^{i\psi(x, \xi_\alpha)} T_\alpha f(x)\Big| +\max_\alpha|T_\alpha f(x)|\tag 2.14
$$
and we write on $B_K$
$$
\sum_{\Omega_\alpha \subset \tilde{\Cal L}} e^{i\psi(x, \xi_\alpha)} T_\alpha f(x) =\phi_{B_K}(x). \Big(\sum_{\Omega_\alpha\subset \tilde{\Cal L}}
|T_\alpha f|^2\Big) ^{1/2}.\tag 2.15
$$
Next, define a function $\phi$ by
$$
\phi|_{B_K}=\phi_{B_K}
$$
for $B_K$ satisfying alternative II.
It follows from (2.14) that the corresponding contribution to the maximal function is bounded by
$$
\Big\Vert\max_{|x_{n+1}|<R}\Big\{ \phi\Big(\sum_\alpha |T_\alpha f|^2 \Big)^{\frac 12}\Big\}\Big\Vert_{L^2[|x'|<R]}.\tag 2.16
$$
Let $\{B_{\gamma, \ell}\}, B_{\gamma, \ell} = B_\gamma \times I_\ell \subset \Bbb R^n\times \Bbb R$ be a partition of $B(0, R)$ in $K$-cubes.
Clearly
$$
(2.16) \leq \Big\{\sum_{\gamma, \ell} \Big(\sum_\alpha |T_\alpha f|^2\Big) \Big|_{B_{\gamma, \ell}} 
 \Big[ \int_{B_\gamma} \max_{x_{n+1} \in I_\ell} |\phi_{\gamma, \ell}|^2
dx'\Big]\Big\}^{1/2}.\tag 2.17
$$
Assume we dispose over a bound
$$
\Big[\nint_{B_\gamma} \max_{x_{n+1}\in I_\ell} |\phi_{\gamma, \ell}|^2 dx'\Big]^{1/2} <A\tag 2.18
$$
where $\nint_{B_\gamma}$ stands for $\frac 1{\mes B_\gamma} \int_{B_\gamma}$.

Then (2.17) is bounded by
$$
\align
&A\Big\{\frac 1K \sum_{\gamma, \ell}\int_{B_{\gamma, \ell}} \Big(\sum_\alpha |T_\alpha f|^2\Big) dx\Big\}^{1/2}\\
&=A \Big(\frac 1K\Big)^{\frac 12} \Big\Vert\Big(\sum_\alpha |T_\alpha f|^2\Big)^{1/2} \Big\Vert_{L^2(|x|<R)}\\
&\leq CA\Big(\frac RK\Big)^{1/2} \Big(\sum_\alpha\Vert f|_{\Omega_\alpha} \Vert^2_2\Big)^{1/2}\\
&= CA\Big(\frac RK\Big)^{1/2} \Vert f\Vert_2.\tag 2.19
\endalign
$$
It remains to establish a bound on $A$ in (2.18).

Again, since in (2.15) the $T_\alpha f$ are viewed as constant on $B_K$, we need a bound
of the form
$$
\Big\Vert \max_{|x_{n+1}|<K}\Big|\sum_\alpha a_\alpha e^{i\psi(x, \xi_\alpha)}\Big|\, \Big\Vert
_{L^2[|x'|<K]}\leq AK^{\frac n2}\Big(\sum |a_\alpha|^2\Big)^{1/2}\tag 2.20
$$
where the points $\{\xi_\alpha\}, |\xi_\alpha|\sim 1$, are $\frac 1K$-separated and in an $(n-1)$-dim affine hyperplane $\Cal L\subset \Bbb R^n$.

Performing a rotation around the $e_{n+1}$-axis, we may assume $\xi_\alpha=c=$ constant for $\xi =(\xi_1, \ldots, \xi_n)\in\Cal L$.
Hence the left side of (2.20) becomes
$$
K^{\frac 12} \Big\Vert\max_{|x_{n+1}|<K}\Big| \sum_\alpha a_\alpha e^{i\big(x_1\xi_{\alpha, 1}+\cdots+x_{n-1} \xi_{\alpha, n-1}
+x_{n+1} (\xi^2_{\alpha, 1} +\cdots+\xi^2_{\alpha, n-1})\big)}\Big| \, \Big\Vert_{L^2_{[|x_1|, \ldots, |x_{n-1}|<K]}}.\tag 2.21
$$
At this point, the dimension is reduced from $n$ to $n-1$.

Let $\{\Omega_\alpha'\}$ be disjoint $\frac 1K$-neighborhoods of the points $\{(\xi_{\alpha, 1}, \ldots, 
\xi_{\alpha, n-1})\}$ in $\Bbb R^{n-1}$.

Let $g$ on $[|\xi'|\sim 1]\subset \Bbb R^{n-1}$ be defined by
$$
g(\xi') =\left\{\aligned &a_\alpha \ |\Omega_\alpha'|^{-1} \ \text { if } \xi'\in \Omega_\alpha'\\
&0 \qquad\qquad\qquad \text { otherwise.}\endaligned
\right.
$$
Hence
$$
\Vert g\Vert_2 \sim K^{\frac {n-1}2} \Big(\sum |a_\alpha|^2\Big)^{1/2}\tag 2.22
$$
and
$$
(2.21)\sim K^{\frac 12} \Big\Vert\max_{|t|<K} \Big|\int g(\xi') e^{i(y.\xi'+t|\xi'|^2)} \Big| \, \Big\Vert_{L^2[|y|<K]}.\tag 2.23
$$
Denote $\theta_{n-1}$ an exponent for which (2.3) holds in dimension $n-1$.

From (2.22), it follows that
$$
(2.23) \lesssim K^{\frac 12+\theta_{n-1}} K^{\frac {n-1}2}\Big(\sum_\alpha|a_\alpha|^2\Big)^{1/2}
$$
which shows that we may take $A\sim K^{\theta_{n-1}}$.
Hence in (2.19), we obtain the estimate
$$
CR^{\frac 12} K^{\theta_{n-1}-\frac 12}\Vert f\Vert_2.\tag 2.24
$$
Together with the Case I contribution (2.13), we obtain the bound
$$
\big\Vert\sup_{0<t<R} |e^{it\Delta} f|\big\Vert_{L^2_{(B_R)}} \ll [R^{\frac n{2(n+1)}+\ve} K^C+R^{\frac 12} K^{\theta_{n-1}-\frac 12}] \Vert
f\Vert_2
\tag 2.25
$$
where $K$ is a parameter.
Recall that $\theta_{n-1} <\frac 12$. 
Hence an appropriate choice of $K=R^\delta$ permits to obtain (2.3) for some $\theta =\theta_n< \frac 12$.
\bigskip

\noindent
{\bf 3. A quantitative estimate}

The proof of the second part of Theorem 1 depends on the more refined analysis from \S4 in \cite {B-G} and its higher dimensional version.

Before stating the relevant result from \cite {B-G}, we fix some terminology.
Let
$$
S=\{(\xi, |\xi|^2), |\xi|\sim 1\} \subset\Bbb R^{n+1}.
$$
For fixed $\delta>0$, consider a partition of $\{\xi:|\xi|\sim 1\}\subset\Bbb R^n$ in cells $Q$ of size $\delta$ and denote
$$
\tau=\{(\xi, |\xi|^2); \xi \in Q\}\tag 3.1
$$
the corresponding partition of $S$ in $\delta$-caps.
Thus the convex hull $\text{conv} (\tau)$ is a $(\underbrace{\delta\times\cdots\times \delta}_n\times \delta^2)$-box tangent to $S$ and we denote
$\overset\circ \to \tau \subset \Bbb R^{n+1}$ the polar box of $\text{conv}(\tau)$ shifted with $0$ as center.
Thus $\overset\circ\to\tau$ is a $\big(\frac 1\delta\times\cdots \times \frac 1\delta\times \frac 1{\delta^2})$-box.
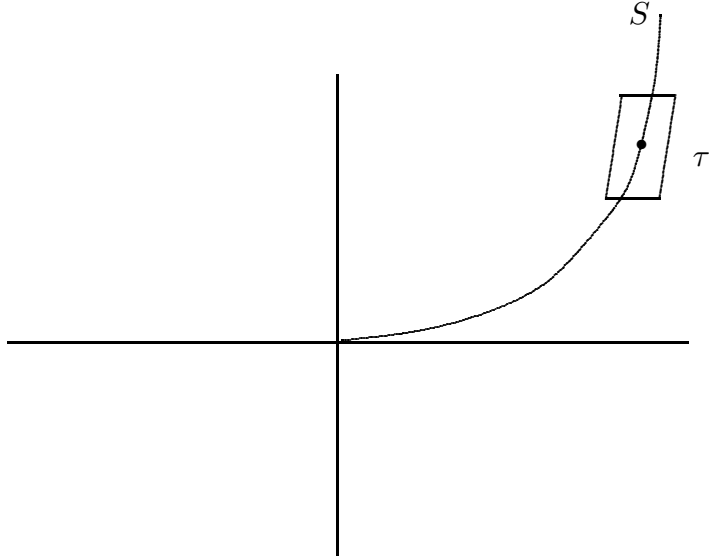
\input pictex
\font\thinlinefont=cmr5
$$
\hbox{\beginpicture
\setcoordinatesystem units <.80000cm,.80000cm>
% Fig POLYLINE object
\linethickness= 0.500pt
\setplotsymbol ({\thinlinefont .})
{\putrule from  7.906 25.495 to  7.906 18.828
}%
% Fig POLYLINE object
\linethickness= 0.500pt
\setplotsymbol ({\thinlinefont .})
{\plot 12.366 23.779 12.588 25.207 /
}%
% Fig POLYLINE object
\linethickness= 0.500pt
\setplotsymbol ({\thinlinefont .})
{\putrule from 11.620 23.781 to 12.383 23.781
}%
% Fig POLYLINE object
\linethickness= 0.500pt
\setplotsymbol ({\thinlinefont .})
{\plot 11.625 23.774 11.847 25.203 /
}%
% Fig POLYLINE object
\linethickness= 0.500pt
\setplotsymbol ({\thinlinefont .})
{\putrule from 11.811 25.210 to 12.573 25.210
}%
% Fig POLYLINE object
\linethickness= 0.500pt
\setplotsymbol ({\thinlinefont .})
{\putrule from  3.334 21.780 to 12.764 21.780
}%
% Fig POLYLINE object
\linethickness= 0.500pt
\setplotsymbol ({\thinlinefont .})
{\putrule from  7.969 21.812 to  7.976 21.812
\plot  7.976 21.812  7.993 21.814 /
\plot  7.993 21.814  8.020 21.816 /
\plot  8.020 21.816  8.060 21.821 /
\plot  8.060 21.821  8.111 21.827 /
\plot  8.111 21.827  8.172 21.833 /
\plot  8.172 21.833  8.240 21.842 /
\plot  8.240 21.842  8.310 21.848 /
\plot  8.310 21.848  8.380 21.857 /
\plot  8.380 21.857  8.450 21.865 /
\plot  8.450 21.865  8.513 21.872 /
\plot  8.513 21.872  8.575 21.878 /
\plot  8.575 21.878  8.630 21.886 /
\plot  8.630 21.886  8.680 21.893 /
\plot  8.680 21.893  8.729 21.899 /
\plot  8.729 21.899  8.776 21.905 /
\plot  8.776 21.905  8.818 21.910 /
\plot  8.776 21.905  8.818 21.910 /
\plot  8.818 21.910  8.860 21.916 /
\plot  8.860 21.916  8.901 21.922 /
\plot  8.901 21.922  8.941 21.931 /
\plot  8.941 21.931  8.983 21.937 /
\plot  8.983 21.937  9.023 21.943 /
\plot  9.023 21.943  9.066 21.952 /
\plot  9.066 21.952  9.108 21.960 /
\plot  9.108 21.960  9.152 21.969 /
\plot  9.152 21.969  9.197 21.979 /
\plot  9.197 21.979  9.243 21.988 /
\plot  9.243 21.988  9.288 21.999 /
\plot  9.288 21.999  9.335 22.009 /
\plot  9.335 22.009  9.381 22.020 /
\plot  9.381 22.020  9.426 22.032 /
\plot  9.426 22.032  9.472 22.043 /
\plot  9.472 22.043  9.517 22.056 /
\plot  9.517 22.056  9.561 22.068 /
\plot  9.561 22.068  9.603 22.079 /
\plot  9.603 22.079  9.646 22.092 /
\plot  9.646 22.092  9.686 22.104 /
\plot  9.686 22.104  9.728 22.117 /
\plot  9.728 22.117  9.768 22.130 /
\plot  9.768 22.130  9.809 22.142 /
\plot  9.809 22.142  9.851 22.157 /
\plot  9.851 22.157  9.893 22.172 /
\plot  9.893 22.172  9.938 22.187 /
\plot  9.938 22.187  9.982 22.204 /
\plot  9.982 22.204 10.027 22.221 /
\plot 10.027 22.221 10.073 22.238 /
\plot 10.073 22.238 10.122 22.257 /
\plot 10.122 22.257 10.168 22.276 /
\plot 10.168 22.276 10.215 22.295 /
\plot 10.215 22.295 10.262 22.316 /
\plot 10.262 22.316 10.308 22.335 /
\plot 10.308 22.335 10.353 22.356 /
\plot 10.353 22.356 10.397 22.375 /
\plot 10.397 22.375 10.437 22.396 /
\plot 10.437 22.396 10.478 22.416 /
\plot 10.478 22.416 10.516 22.437 /
\plot 10.516 22.437 10.552 22.456 /
\plot 10.552 22.456 10.588 22.475 /
\plot 10.588 22.475 10.621 22.494 /
\plot 10.621 22.494 10.655 22.517 /
\plot 10.655 22.517 10.691 22.540 /
\plot 10.691 22.540 10.725 22.564 /
\plot 10.725 22.564 10.761 22.587 /
\plot 10.761 22.587 10.795 22.612 /
\plot 10.795 22.612 10.829 22.640 /
\plot 10.829 22.640 10.863 22.667 /
\plot 10.863 22.667 10.897 22.697 /
\plot 10.897 22.697 10.930 22.727 /
\plot 10.930 22.727 10.964 22.756 /
\plot 10.964 22.756 10.998 22.788 /
\plot 10.998 22.788 11.030 22.820 /
\plot 11.030 22.820 11.060 22.852 /
\plot 11.060 22.852 11.091 22.881 /
\plot 11.091 22.881 11.121 22.913 /
\plot 11.121 22.913 11.151 22.945 /
\plot 11.151 22.945 11.178 22.976 /
\plot 11.178 22.976 11.208 23.008 /
\plot 11.208 23.008 11.233 23.038 /
\plot 11.233 23.038 11.261 23.067 /
\plot 11.261 23.067 11.288 23.099 /
\plot 11.288 23.099 11.318 23.133 /
\plot 11.318 23.133 11.347 23.167 /
\plot 11.347 23.167 11.377 23.201 /
\plot 11.377 23.201 11.409 23.239 /
\plot 11.409 23.239 11.441 23.275 /
\plot 11.441 23.275 11.472 23.313 /
\plot 11.472 23.313 11.502 23.349 /
\plot 11.502 23.349 11.534 23.387 /
\plot 11.534 23.387 11.563 23.425 /
\plot 11.563 23.425 11.593 23.461 /
\plot 11.593 23.461 11.623 23.497 /
\plot 11.623 23.497 11.648 23.531 /
\plot 11.648 23.531 11.676 23.563 /
\plot 11.676 23.563 11.699 23.594 /
\plot 11.699 23.594 11.722 23.624 /
\plot 11.722 23.624 11.743 23.654 /
\plot 11.743 23.654 11.764 23.679 /
\plot 11.764 23.679 11.788 23.713 /
\plot 11.788 23.713 11.809 23.745 /
\plot 11.809 23.745 11.832 23.779 /
\plot 11.832 23.779 11.853 23.810 /
\plot 11.853 23.810 11.872 23.844 /
\plot 11.872 23.844 11.891 23.880 /
\plot 11.891 23.880 11.913 23.916 /
\plot 11.913 23.916 11.930 23.954 /
\plot 11.930 23.954 11.949 23.992 /
\plot 11.949 23.992 11.966 24.033 /
\plot 11.966 24.033 11.980 24.073 /
\plot 11.980 24.073 11.997 24.115 /
\plot 11.997 24.115 12.012 24.158 /
\plot 12.012 24.158 12.025 24.204 /
\plot 12.025 24.204 12.040 24.251 /
\plot 12.040 24.251 12.054 24.299 /
\plot 12.054 24.299 12.065 24.337 /
\plot 12.065 24.337 12.076 24.378 /
\plot 12.076 24.378 12.088 24.420 /
\plot 12.088 24.420 12.099 24.464 /
\plot 12.099 24.464 12.112 24.513 /
\plot 12.112 24.513 12.124 24.562 /
\plot 12.124 24.562 12.137 24.615 /
\plot 12.137 24.615 12.150 24.668 /
\plot 12.150 24.668 12.162 24.723 /
\plot 12.162 24.723 12.175 24.780 /
\plot 12.175 24.780 12.188 24.839 /
\plot 12.188 24.839 12.200 24.896 /
\plot 12.200 24.896 12.213 24.956 /
\plot 12.213 24.956 12.226 25.013 /
\plot 12.226 25.013 12.239 25.070 /
\plot 12.239 25.070 12.249 25.125 /
\plot 12.249 25.125 12.260 25.180 /
\plot 12.260 25.180 12.270 25.233 /
\plot 12.270 25.233 12.279 25.284 /
\plot 12.279 25.284 12.287 25.332 /
\plot 12.287 25.332 12.296 25.381 /
\plot 12.296 25.381 12.304 25.425 /
\plot 12.304 25.425 12.311 25.476 /
\plot 12.311 25.476 12.317 25.525 /
\plot 12.317 25.525 12.325 25.574 /
\plot 12.325 25.574 12.330 25.622 /
\plot 12.330 25.622 12.336 25.671 /
\plot 12.336 25.671 12.342 25.724 /
\plot 12.342 25.724 12.347 25.777 /
\plot 12.347 25.777 12.351 25.836 /
\plot 12.351 25.836 12.355 25.895 /
\plot 12.355 25.895 12.361 25.957 /
\plot 12.361 25.957 12.366 26.020 /
\plot 12.366 26.020 12.370 26.084 /
\plot 12.370 26.084 12.372 26.143 /
\plot 12.372 26.143 12.376 26.198 /
\plot 12.376 26.198 12.378 26.242 /
\plot 12.378 26.242 12.380 26.278 /
\plot 12.380 26.278 12.383 26.302 /
\putrule from 12.383 26.302 to 12.383 26.314
\putrule from 12.383 26.314 to 12.383 26.321
}%
%
% Fig TEXT object
%
\put{$_\bullet$} [lB] at 12.029 24.511
%
% Fig TEXT object
%
\put{$S$} [lB] at 11.938 26.194
%
% Fig TEXT object
%
\put{$\tau$} [lB] at 12.827 24.225
\linethickness=0pt
\putrectangle corners at  3.308 26.480 and 13.411 18.802
\endpicture}
$$

Denote also
$$
f_\tau=f|_Q
$$
with $Q$ and $\tau$ related by (3.1).

Let the  operator $T$ be as in previous section.

Next recall (3.4) - (3.8) from \S4 in \cite {B-G}, providing an estimate for $Tf$ on $B_R$ in $3D$ (i.e. $n=2)$. Thus
$$
\align
&|Tf| \ll\\
&R^\ve \max_{\frac 1{\sqrt T}<\delta<1} \max_{\Cal E_\delta} \Big[\sum_{\tau\in\Cal E_\delta} (\phi_\tau |Tf_{\tau_1}|^{1/3} |Tf_{\tau_2}|^{1/3}
|Tf_{\tau_3}|^{1/3})^2\Big]^{1/2}\\
&R^\ve \max_{\Cal E_{\frac 1{\sqrt R}}} \Big[\sum_{\tau\in\Cal E} (\phi_\tau |Tf_{\tau})^2 \Big]^{1/2}\tag 3.2
\endalign
$$
where
\roster
\item"{(3.3)}" $\Cal E_\delta$ consists of at most $\frac 1\delta$ disjoint $\delta$-caps.

\item "{(2.4)}" $\tau_1, \tau_2, \tau_3\subset\tau$ are 3-transversal $\frac \delta K$-caps ($K$ is  a large constant).

\item "{(3.5)}" For each $\tau, \phi_\tau\geq 0$ is a function on $\Bbb R^n$ satisfying
$$
\nint_{B}\phi^4_\tau \ll R^\ve
$$
for all $B$ taken in a tiling of $\Bbb R^{n+1}$ with translates of $\overset\circ\to \tau$.
\endroster

Of course (3.2) certainly implies
$$
\align
|Tf|&\ll R^\ve \sum_{\Sb \delta\text { dyadic}\\ \frac 1{\sqrt R}<\delta<1\endSb}
\Big[\sum_{\tau \, \delta-\text{cap}}(\phi_\tau |Tf_{\tau_1}| ^{1/3} |Tf_{\tau_2}|^{1/3} |Tf_{\tau_3}|^{1/3})^2\Big]^{1/2}\tag 3.6\\
&+ R^\ve \Big[\sum_{\tau\frac 1{\sqrt R}-\text{cap}} (\phi_\tau |Tf_\tau|)^2\Big]^{1/2}
\tag 3.7
\endalign
$$
where $\sum\limits_{\tau \ \delta-\text{cap}}$ refers to summation over a partition of $S$ in $\delta$-caps $\tau$.

Formula (3.2) has a higher dimensional version, deduced in a similar way from \S3 in \cite{B-G}.
In particular, we get for general $n\geq 2$
$$
\align
|Tf|&\ll\\
&R^\ve\sum_\delta\Big[\sum_{\tau \, \delta-\text{cap}} \Big(\phi_\tau \prod^{n+1}_{j=1} |Tf_{\tau_j}|^{\frac 1{n+1}}\Big)^2 
\Big]^{\frac 12}\tag 3.8\\
&+ R^\ve\Big[\sum_{\tau\frac 1{\sqrt R}-\text{cap}} (\phi_\tau|Tf_\tau|)^2\Big]^{\frac 12}\tag 3.9
\endalign
$$
where now $\tau_1, \ldots, \tau_{n+1} \subset \tau$ are $(n+1)$-transversal $\frac \delta K$-caps and $\phi_\tau\geq 0$ on $\Bbb R^{n+1}$ satisfies
$$
\nint_B \phi_\tau^q \ll R^\ve \text { with } q=\frac {2n}{n-1}\tag 3.10
$$
for all $B$ in a tiling of $\Bbb R^{n+1}$ from $\overset\circ\to \tau$-translates.

It is now clear from (3.8), (3.9) and convexity that in order to bound
$$
\Vert Tf\Vert_{L^2_{[|x'|<R]} L^\infty_{[|x_{n+1}|<R]}}
$$
it suffices to estimate
$$
\Big\Vert\phi_\tau\Big(\prod_{j=1}^{n+1} |Tf_{\tau_j}|^{\frac 1{n+1}}\Big)\Big\Vert_{L^2_{[|x'|<R]} L^\infty_{[|x_{n+1}|<R]}}\tag 3.11
$$
with $\tau_1, \ldots, \tau_{n+1}\subset \tau$ as above, and also
$$
\Vert\phi_\tau|Tf_\tau| \, \Vert_{L^2_{[|x'|<R]}L^\infty_{[|x_{n+1}|<R]}}\tag 3.12
$$
with $\tau$ or $\frac 1{\sqrt R}$-cap.

Denote $\phi=\phi_\tau$ and $F=\prod^{n+1}_{j=1} |Tf_{\tau_j}|^{\frac 1{n+1}}$.

Let $\{B_{jk}\}$ be a tiling of $B(0, R)\subset \Bbb R^{n+1}$ with $\pi_{x'}(B_{jk})=I_j\subset\Bbb R^n$ a partition of $B(0, R)\subset \Bbb R^n$
in $\big(\underbrace{\tfrac 1\delta\times\cdots\times \tfrac1\delta}_{n-1} \times \frac 1{\delta^2})$-boxes and $B_{jk} \sim$ translate of
$\overset\circ\to\tau$
\font\thinlinefont=cmr5
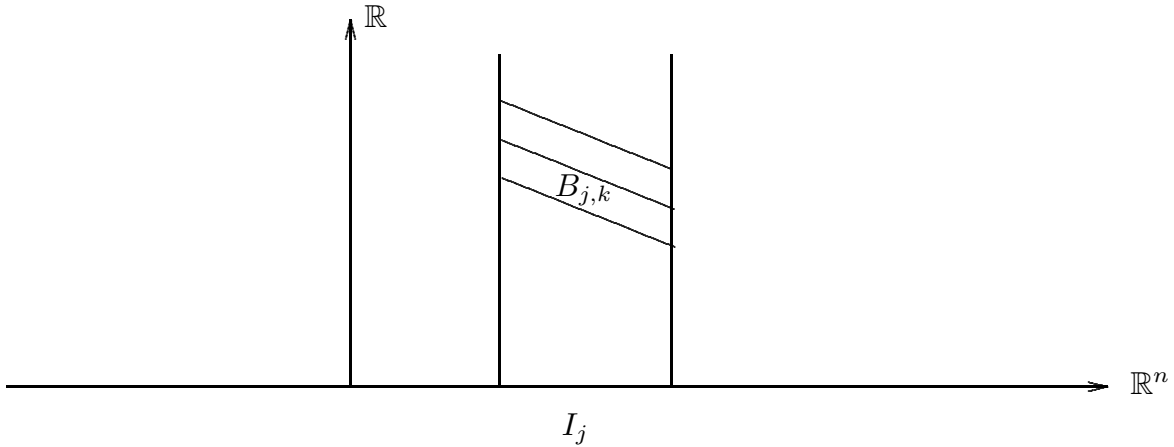
$$
\hbox{\beginpicture
\setcoordinatesystem units <0.80000cm,.80000cm>
% Fig POLYLINE object
\linethickness= 0.500pt
\setplotsymbol ({\thinlinefont .})
{\putrule from  6.032 20.955 to  6.032 26.035
% arrow head
\plot  6.096 25.781  6.032 26.035  5.969 25.781 /
}%
% Fig POLYLINE object
\linethickness= 0.500pt
\setplotsymbol ({\thinlinefont .})
{\putrule from  8.096 20.955 to  8.096 25.559
}%
% Fig POLYLINE object
\linethickness= 0.500pt
\setplotsymbol ({\thinlinefont .})
{\putrule from 10.478 20.955 to 10.478 25.559
}%
% Fig POLYLINE object
\linethickness= 0.500pt
\setplotsymbol ({\thinlinefont .})
{\putrule from  1.270 20.955 to 16.510 20.955
% arrow head
\plot 16.256 20.892 16.510 20.955 16.256 21.018 /
}%
% Fig POLYLINE object
\linethickness= 0.500pt
\setplotsymbol ({\thinlinefont .})
{\plot  8.096 24.924 10.478 23.971 /
}%
% Fig POLYLINE object
\linethickness= 0.500pt
\setplotsymbol ({\thinlinefont .})
{\plot  8.130 24.371 10.511 23.419 /
}%
% Fig POLYLINE object
\linethickness= 0.500pt
\setplotsymbol ({\thinlinefont .})
{\plot  8.136 23.846 10.518 22.894 /
}%
% Fig TEXT object
\put{$I_j$} [lB] at  8.954 20.288
% Fig TEXT object
\put{$B_{j,k}$} [lB] at  8.858 23.590
% Fig TEXT object
\put{$\Bbb R$} [lB] at  6.223 25.940
% Fig TEXT object
\put{$\Bbb R^n$} [lB] at 16.828 20.860
\linethickness=0pt
\putrectangle corners at  1.245 26.226 and 17.217 20.185
\endpicture}
$$
~

Let
$$
p=\frac {2(n+1)}{n} \ \text { and } \  q=\frac {2n}{n-1}.
$$
By H\"older's inequality
$$
\align
(3.11) &= \Big[\sum_j \Vert\phi F\Vert^2_{L^2_{I_j}L^\infty_{[|x_{n+1}|<R]}}\Big]^{\frac 12}\\
& \leq \delta^{-\frac 12}\Big[\sum_j \Vert\phi F\Vert^2_{L^p_{I_j}L^\infty_{x_{n+1}}}\Big]^{\frac 12}.\tag 3.13
\endalign
$$
Writing for fixed $x'\in I_j$
$$
\sup_{x_{n+1}} |\phi F|(x', x_{n+1})\leq \Big(\sum_k\Vert\phi F\Vert^p_{L^\infty(B_{j, k}(x'))}\Big)^{\frac 1p}
$$
it follows that
$$
(3.13) \leq \delta^{-\frac 12}\Big[\sum_j \Big(\sum_k\Vert\phi F\Vert^p_{L^p_{x'} L^\infty _{x_{n+1}}(B_{j, k})}\Big)^{\frac 2p}\Big]^{\frac
12}.\tag 3.14
$$
Since $\supp \widehat {Tf_{\tau_j}} \subset\tau$, we may view $F$ as essentially constant on each $B_{jk}$.
From (3.10)
$$
\align
\Vert\phi\Vert_{L^p_{x'} L^\infty_{x_{n+1}}(B_{jk})} &\leq \Vert\phi\Vert_{L^p_{x'} L^q_{x_{n+1}(B_{jk})}}\\
&\leq |I_j|^{\frac 1p-\frac 1q} \Vert\phi\Vert_{L^q(B_{jk})}\\
&\ll R^\ve |I_j|^{\frac 1p-\frac 1q} |B_{jk}|^{\frac 1q}\\
&\ll R^\ve \Big(\frac 1\delta\Big)^{\frac {n}2+\frac {n-1}{2n}}
\endalign
$$
and hence
$$
\spreadlines{5pt}
\align
\Vert\phi F\Vert_{L^p_{x'} L^\infty_{x_{n+1}(B_{jk})}}&\ll R^\ve\Big(\frac 1\delta\Big) ^{\frac n2+\frac{n-1}{2n}}|I_j|^{-\frac 12} 
\delta^{\frac 1p}\Vert F\Vert_{L^2_{x'} L^p_{x_{n+1}}(B_{jk})}\\
&\ll R^\ve \delta^{\frac 12+ \frac 1{2n} -\frac 1{2(n+1)}} \Vert F\Vert_{L^2_{x'}L^p_{x_{n+1}}(B_{jk})}.\tag 3.15
\endalign
$$
Substitution of (3.15) in (3.14) gives
$$
\spreadlines{5pt}
\align
&R^\ve \delta^{\frac 1{2n}-\frac 1{2(n+1)}} \Big[\sum_j \Big(\sum_k \Vert F\Vert^p_{L^2_{x'} L^p_{x_{n+1}}(B_{jk})} \Big)^{\frac 2p}\Big]^{\frac
12}\\
&\leq R^\ve \delta^{\frac 1{2n}-\frac 1{2(n+1)}} \Vert F\Vert_{L^2_{[|x'|<R]} L^p_{[|x_{n+1}|< R]}}.\tag 3.16
\endalign
$$
Next, we evaluate $\Vert F\Vert_{L^2_{x'} L^p_{x_{n+1}}}$.
Recall the definition of $F$.

Let $\xi_0$ be the center of $Q=\pi_\xi(\tau)\subset\Bbb R^n$ and write $\xi-\xi_0=\delta\zeta, |\zeta|=O(1)$, for $\xi\in Q$.
Hence, rescaling
$$
\spreadlines{5pt}
\align
|(Tf_{\tau_j})(x)|&= \delta^n \Big|\int e^{i((\delta x'+ 2\delta\xi_0 x_{n+1}).\zeta +\delta^2|\zeta|^2 x_{n+1})}  
f_{\tau_j} (\xi_0+\delta\zeta) d\zeta\Big|\\
&=\delta^n\Big|\int e^{i((y'+2\delta^{-1} y_{n+1} \xi_0). \zeta+|\zeta|^2 y_{n+1})} f_{\tau_j} (\xi_0+\delta\zeta) d\zeta\Big|\\
&=\delta^{\frac n2} |Tg_{\tilde \tau_j}| (y'+2\delta^{-1} y_{n+1} \xi_0, y_{n+1})\tag 3.17
\endalign
$$
where $(y_1, \ldots, y_{n+1})= (\delta x_1, \ldots, \delta x_n, \delta^2 x_{n+1})$ satisfies
$$
|y'|<\delta R, |y_{n+1}|< \delta^2 R
$$
and $\Vert g_{\tilde\tau_j}\Vert_2=1$.

Denote
$$
G=\prod^{n+1}_{1} |Tg_{\tilde\tau_j} |^{\frac 1{n+1}}
$$
where $\tilde\tau_1, \ldots, \tilde\tau_{n+1}$ are transversal $O(1)$-caps.

Thus
$$
\Vert F\Vert_{L^2_{x'} L^p_{x_{n+1}}} =\delta^{-\frac 2p} \Vert G(y'+2\delta^{-1} y_{n+1} \xi_0, y_{n+1})\Vert _{L^2_{y'} 
L^p_{y_{n+1}}}.\tag 3.18
$$
Assume (as we may by rotation in $[e_1, \ldots, e_n]$-space) that $\xi_0// e_1$.

Estimate using again H\"older's inequality
$$
\Vert G(y_1+2\delta^{-1} y_{n+1} \xi_0, y_2, \ldots, y_n, y_{n+1})\Vert_{L^2_{y'} L^p_{y_{n+1}}} \leq (\delta R)^{\frac 12-\frac 1p}\Vert
G\Vert_{L^2_{y_2, \ldots, y_n}L^p_{y_1, y_{n+1}}}.
$$
This gives by (3.18)
$$
(3.16) < \delta^{\frac 1{2n}-\frac n{n+1}} R^{\frac 1{2(n+1)}+\ve} \Vert G\Vert_{L^2_{y_2, \ldots, y_n} L^p_{y_1, y_{n+1}}}\tag 3.19
$$
where $|y_1|, \ldots, |y_n|<\delta R, |y_{n+1}| <\delta^2R$.

Partition $[|y'|< \delta R]$ in cubes $\Omega_s$ of size $\delta^2 R$.
Clearly
$$
\Vert G\Vert_{L^2_{y_2, \ldots, y_n}L^p_{y_1, y_{n+1}}} \leq (\delta^2 R)^{(n-1)(\frac 12-\frac 1p)} \Big[\sum_s
\Vert G\Vert^2_{L^p_{y' \in \Omega_s}L^p_{y_{n+1}}}\Big]^{\frac 12}.\tag 3.20
$$
Recalling the definition of $G$ and using again the \cite{B-C-T} multi-linear restriction inequality (together with a standard localization
argument in $y'$-space), it follows that
$$
\Big[\sum_s\Vert G\Vert^2_{L^p_{y'\in\Omega_s} L^p_{y_{n+1}}} \Big]^{\frac 12}\ll R^\ve.
$$
Hence
$$
(3.16)< \delta^{\frac 1{2n}-\frac 1{n+1}} R^{\frac n{2(n+1)}+\ve} \tag 3.21
$$
and since $\delta>\frac 1{\sqrt R}$, this implies that
$$
(3.11)\ll R^{\frac 12 -\frac 1{4n}+\ve}.\tag 3.22
$$
Next, consider (3.12) for which we repeat the preceding with $\delta=\frac 1{\sqrt R}$ up to (3.20), with $G=|Tg_{
\tilde \tau}|, \Vert g_{\tilde\tau}\Vert_2 =1$.
The cubes $\Omega_s$ and $y_{n+1}$ are size $O(1)$.
Hence the $RH$-side of (3.20) becomes
$$
\Vert G\Vert_{L^2_{[|y'|<\sqrt R]} L^2_{[|y_{n+1}| < 0(1)]}}\lesssim \Vert g_{\tilde\tau}\Vert_2 =1
$$
and (3.12) is also bounded by (3.22).

This concludes the proof of Theorem 1.

\bigskip

\noindent
{\bf 4. A construction in higher dimension}

We prove Theorem 2.

Let $S=\{(\xi, \frac 1R|\xi|^2); \xi\in\Bbb R^n, |\xi|\sim R\}$ with $R$ large and let $H\subset\Bbb R^{n+1}$ be the hyperplane
$[z_1-z_{n+1}]$.
Hence $\pi_\xi(S\cap H)\subset\Bbb R ^n$ is the $(n-1)$-sphere
$$
\Big(\xi_1-\frac R2\Big)^2+\xi^2_2+\cdots +\xi _n^2 =\frac {R^2} 4.\tag 4.1
$$
Consider the lattice points
$$
\Cal E=\Bbb Z^n\cap R^{\frac 1n} S^{(n-1)}\tag 4.2
$$
where $S^{(n-1)}$ denotes the unit sphere in $\Bbb R^n$.
If $U \in O(n)$ is an arbitrary orthogonal transformation, we have
$$
\frac 12 R^{\frac {n-1}n} U(\Cal E)\subset\frac R2 S^{(n-1)}
$$
and therefore
$$
\Cal E_1 =\frac R2 e_1+ \frac 12 R^{\frac {n-1}n} U(\Cal E) \subset \pi_\xi(S\cap H).\tag 4.3
$$
Hence $(\xi, \frac 1R|\xi|^2)\in S\cap H$ for $\xi\in\Cal E_1$.
Note also that the points of $\Cal E_1$ are $\sim R^{\frac {n-1}n}$-separated.

Define next
$$
\widehat f =\frac 1{|\Cal E_1|^{\frac 12}} \sum_{\xi\in\Cal E_1} 1_{B(\xi, \rho)}\tag 4.4
$$
with $\rho$ a sufficiently small constant.
Hence $\Vert f\Vert_2=O (1)$ and $\supp\widehat f \subset B_R$.

Let $|x|, |t|< C$.
By construction and taking $\rho$ small enough
$$
(e^{i\frac tR\Delta}f)(x) =\frac 1{|\Cal E|^{\frac 12}} \sum_{\xi\in\Cal E_1} \Big[\int_{B(0, \rho)} e^{i(x. \zeta +\frac tR|\xi+\zeta|^2)}
d\zeta\Big] e^{ix.\xi}
$$
and
$$
|(e^{i\frac tR\Delta} f)(x)| \sim|\Cal E|^{\frac 12}\tag 4.5
$$
provided $(x, t)$ satisfies moreover
$$
x.\xi+\frac tR|\xi|^2 \in\frac R 2 x_1+ 2\pi\Bbb Z+ B\Big(0, \frac 1{10}\Big)\text { for all } \xi\in\Cal E_1.\tag 4.6
$$
Denote $\nu=\frac 1{\sqrt 2}(e_1-e_{n+1})$ the normal of $H$.
Since $\big(\xi, \frac 1R|\xi|^2)\in H$, it follows that if $(x, t)$ satisfies (4.6), then so does $(x, t)+\Bbb R \nu$.

Also, by (4.3)
$$
x.\xi=\frac R2 x_1 +\frac 12 R^{\frac {n-1}n} U\xi'.x
$$
for some $\xi'\in\Cal E\subset \Bbb Z^n$.
Hence $(x, 0)$ will satisfy (4.6) if $x$ belongs to the dual lattice $\Cal L^* =4\pi R^{-\frac {n-1}n} U(\Bbb Z^n)$ of $\Cal L=
\frac 1{4\pi} R^{-\frac {n-1}n} U(\Bbb Z^n)$.
It follows that (4.6) is valid for
$$
(x, t) \in (\Cal L^*\times\{0\})+\Bbb R\nu+B\Big(0, \frac 1{100 R}\Big).
\tag 4.7
$$
Our aim is to choose $U$ as to ensure that
$$
\pi_x\big((\Cal L^*\times\{0\})+[-C, C]\nu\big) =\Cal L^* +\Big[-\frac C{\sqrt 2} , \frac C{\sqrt 2}\Big] e_1\tag 4.8
$$
is $\frac 1{100 R}$-dense in $B_1\subset \Bbb R^n$.
Recalling the definition of $\Cal L^*$, the issue is to obtain a unit vector $\theta\in\Bbb R^n$ so that
$$
R^{-\frac {n-1}n}\Bbb Z^n +[-C, C]\theta\tag 4.10
$$
is $10^{-4} R^{-1}$ = dense.
Equivalently, given $x\in\Bbb R^n$, there is $\lambda\in\Bbb R$, \break $|\lambda|< O(R^{\frac {n-1}{n}})$ such that
$$
\max_{1\leq j\leq n} \Vert x_j-\lambda\theta_j\Vert < 10^{-4} R^{-\frac 1n}.\tag 4.11
$$
We claim that there is $\theta\in\Bbb R^n, |\theta|=1$ satisfying this property (in fact, the typical $\theta$ will do).
For completeness sake, we will include a proof of this (standard) fact at the end of this section.

The conclusion of the above is that
$$
\big\{ x; \max_{|t|<C}\big|(e^{i\frac tR\Delta} f)(x) \big| \gtrsim |\Cal E|^{\frac 12}\big\} \supset B_1.\tag 4.12
$$
Remains to obtain a lower bound on $|\Cal E|$.
Assume $R^{\frac 1n}\in\Bbb Z$.
If $n\geq 4$, it is well known that for $L^2\in\Bbb Z_+, L\to\infty$
$$
|\Bbb Z^n\cap LS^{(n-1)}|\gtrsim L^{n-2}.\tag 4.13
$$
Thus we may ensure that
$$
|\Cal E|\gtrsim R^{\frac {n-2}n}\tag 4.14
$$
and the implication of (4.12) is that the almost convergence property may only hold for $s\geq \frac {n-2}{2n}$ for $n\geq 4$.
This establishes Theorem 2.

\medskip

\noindent
{\bf Proof of the claim}

Denote $\Bbb T^n=(\Bbb R/\Bbb Z)^n$ the $n$-dim torus.
Let $0\leq \vp\leq 1$ on $\Bbb T^n$ be a smooth function satisfying
$$
\left\{
\aligned
\vp(y) =1 &\text { if } \Vert y\Vert <\frac{R^{-\frac 1n}}{2.10^4}\\
\vp (y) =0 &\text { if } \Vert y\Vert >\frac {R^{-\frac 1n}}{10^4}
\endaligned
\right.
\tag 4.15
$$
and with Fourier expansion
$$
\vp(y) =\sum_{k\in\Bbb Z^n} \hat\vp (k) e(ky)
$$
satisfying
$$
|\hat\vp(k)|\leq C_1 R^{-1} (1+R^{-\frac 1n} |k|)^{-n}.\tag 4.16
$$
To prove (4.11), we show that given any $a\in\Bbb T^n$, there is $\lambda\in\Bbb R, |\lambda|<O\big(R^{\frac {n-1}{n}})
$ s.t
$$
0<\vp(a+\lambda\theta)=\hat\vp(0)+\sum_{k\not= 0}\hat\vp (k) e(k.a) e(\lambda k.\theta).\tag 4.17
$$
Proceed by averaging over $\lambda$ as above.  The contribution of the last term in (4.17) may clearly be bounded by
$$
\sum_{k\not=0} |\hat\vp(k)|\Big|\int e(\lambda k.\theta)\eta (\lambda) d\lambda\Big| \overset {(4.16)}\to < C_1 R^{-1} \sum_{k\not= 0}
(1+R^{-\frac 1n}|k|)^{-n} |\hat\eta (k.\theta)|\tag 4.18
$$
where we let $\eta(u) =R^{-\frac {n-1}n}\eta_1(R^{-\frac {n-1}n} u)$ with $0\leq\eta_1\leq 1, \int_{\Bbb R}\eta_1 =1$ a smooth
compactly supported bumpfunction.
We may ensure that for $|u|\to\infty$
$$
|\hat\eta_1(u)|\leq (1+C_3|u|)^{-2}\tag 4.19
$$
so that
$$
(4.18)< C_1 R^{-1} \sum_{k\not=0}(1+R^{-\frac 1n}|k|)^{-n}(1+C_3R^{\frac {n-1}n} |k. \theta|)^{-2}.\tag 4.20
$$
It will therefore suffice to obtain $\theta$ satisfying $(4.20)<\frac 12 \hat\vp(0)$, i.e.
$$
\sum_{k\not= 0} (1+R^{-\frac 1n}|k|)^{-n} (1+C_3 R^{\frac {n-1}n}|k.\theta|)^{-2} <\frac {C(n)}{C_1}\tag 4.21
$$
with $c(n)$ a constant depending on $n$.
This last statement is easily established by integration over $\theta$ and taking $C_3$ large enough.

\bigskip

\noindent
{\bf 5. Maximal functions associated to oscillatory integral operators}

Returning to Theorem 1, its proof extends to bounding certain maximal functions associated to oscillatory integral operators.

More specifically, assume $\vp:\Omega\to\Bbb R$ a smooth function defined on a neighborhood $\Omega$ of $0\in\Bbb R^n$ of the form
$$
\vp(\xi) =\langle A\xi, \xi\rangle +O(|\xi|^3)\tag 5.1
$$
with $A$ positive definite.

Define
$$
( T^*_Rf)(x')=\sup_{|x_{n+1}|\leq R} \Big|\int_\Omega f(\xi) e^{i(x'.\xi +x_{n+1}\vp (\xi))} d\xi|\tag 5.2
$$
with $f\in L^2(\Omega)$, $x'=(x_1, \ldots, x_n), |x'|\leq R$.

We proved that
$$
\Vert T^*_R f\Vert_{L^2_{[|x'|\leq R]}} \ll R^{\frac 12-\frac 1{4n}+\ve} \Vert f\Vert_2\tag 5.3
$$
and may ask for the best constant $B(R)$ s.t.
$$
\Vert T_R^*f\Vert_{L^2_{[x'|\leq R]}}\leq B(R)\Vert f\Vert_2\tag 5.4
$$
holds in this generality.

Equivalently, let $S$ be a compact, smooth $n$-dim hypersurface in $\Bbb R^{n+1}$ of positive curvature.
Let $\sigma$ be its surface measure.
Fix a point $p\in S$;  let $p+H$ be the tangent hyperplane and $\tau$ the normal at $p$.
Let $\mu\in L^2(S, d\sigma)$, $\Vert\mu\Vert_2=1$ be supported by a neighborhood of $p$ on $S$.
For $R$ large, define $(\hat\mu)^*$ on $B_R\cap H$ as
$$
(\hat\mu)^* (x')=\sup_{|t|<R} |\hat\mu(x'+t\tau)|\tag 5.5
$$
and estimate
$$
\Vert(\hat\mu)^*\Vert _{L^2(H\cap B_R)}.\tag 5.6
$$
Our aim is to describe some examples for which (5.6) is even larger than $R^{\frac 12-\frac 1n}$ as obtained in Theorem 2, hence setting a further limitation to
what one may hope to prove in this general setting.

Consider in the hyperplane $V=[z_1=\ve z_{n+1}]$ a smooth $(n-1)$-dim oval $\Gamma$ with positive curvature, such that $0\in\Gamma$ and $[e_2, \ldots, e_n]$ is the tangent
hyperplane of $\Gamma$ at $0$.
We may clearly produce a hypersurface $S$ in $\Bbb R^{n+1}$ as above, such that $S\cap V=\Gamma$ and moreover the normal $\nu=\frac 1{\sqrt{1+\ve^2}} (-e_1+\ve e_{n+1})$ of
$V$ is tangent to $S$ at each point of $\Gamma$.
Hence, the normal of $S$ at $0$ equals $-\frac 1{\sqrt{1+\ve^2}}(\ve e_1+e_{n+1})$ and for $\ve$ small enough, the point $p\in S$ with normal $\tau=-e_{n+1}$ will be close to
$0$.
Thus $H=[e_2, \ldots, e_n]$
\font\thinlinefont=cmr5
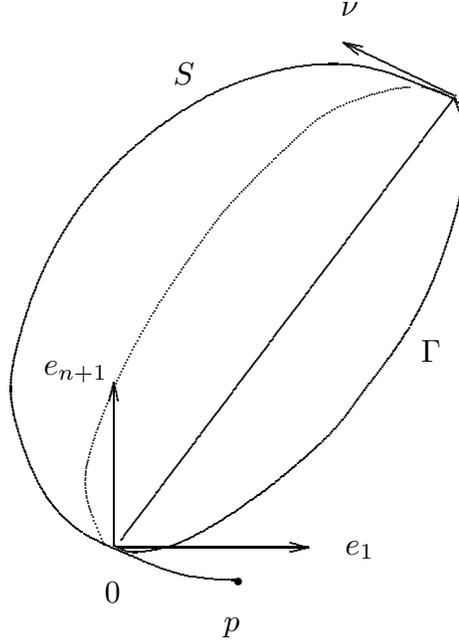
$$
\hbox{\beginpicture
\setcoordinatesystem units <.80000cm,.80000cm>
% Fig POLYLINE object
\linethickness= 0.500pt
\setplotsymbol ({\thinlinefont .})
{
% arrow head
\plot 15.600 22.026 15.344 22.079 15.545 21.912 /
\plot 15.344 22.079 16.868 21.349 /
}%
% Fig POLYLINE object
\linethickness= 0.500pt
\setplotsymbol ({\thinlinefont .})
{\plot 16.895 21.340 12.260 15.181 /
}%
% Fig POLYLINE object
\linethickness= 0.500pt
\setplotsymbol ({\thinlinefont .})
{\putrule from 12.160 15.113 to 12.160 17.367
% arrow head
\plot 12.224 17.113 12.160 17.367 12.097 17.113 /
}%
% Fig POLYLINE object
\linethickness= 0.500pt
\setplotsymbol ({\thinlinefont .})
{\putrule from 12.224 15.081 to 14.859 15.081
% arrow head
\plot 14.605 15.018 14.859 15.081 14.605 15.145 /
}%
% Fig POLYLINE object
\linethickness= 0.500pt
\setplotsymbol ({\thinlinefont .})
\setdots < 0.0953cm>
{\plot 16.954 21.241 16.952 21.241 /
\plot 16.952 21.241 16.946 21.247 /
\plot 16.946 21.247 16.927 21.258 /
\plot 16.927 21.258 16.899 21.277 /
\plot 16.899 21.277 16.863 21.298 /
\plot 16.863 21.298 16.828 21.321 /
\plot 16.828 21.321 16.792 21.340 /
\plot 16.792 21.340 16.760 21.359 /
\plot 16.760 21.359 16.730 21.376 /
\plot 16.730 21.376 16.703 21.389 /
\plot 16.703 21.389 16.677 21.399 /
\plot 16.677 21.399 16.654 21.410 /
\plot 16.654 21.410 16.626 21.419 /
\plot 16.626 21.419 16.599 21.429 /
\plot 16.599 21.429 16.571 21.435 /
\plot 16.571 21.435 16.540 21.444 /
\plot 16.540 21.444 16.506 21.450 /
\plot 16.506 21.450 16.472 21.457 /
\plot 16.472 21.457 16.436 21.461 /
\plot 16.436 21.461 16.400 21.463 /
\plot 16.256 21.463 16.228 21.461 /
\plot 16.228 21.461 16.199 21.459 /
\plot 16.199 21.459 16.169 21.455 /
\plot 16.169 21.455 16.137 21.450 /
\plot 16.137 21.450 16.104 21.444 /
\plot 16.104 21.444 16.070 21.438 /
\plot 16.070 21.438 16.032 21.431 /
\plot 16.032 21.431 15.996 21.425 /
\plot 15.996 21.425 15.958 21.416 /
\plot 15.958 21.416 15.919 21.406 /
\plot 15.919 21.406 15.881 21.397 /
\plot 15.881 21.397 15.843 21.387 /
\plot 15.843 21.387 15.807 21.378 /
\plot 15.807 21.378 15.771 21.368 /
\plot 15.771 21.368 15.735 21.357 /
\plot 15.735 21.357 15.701 21.347 /
\plot 15.701 21.347 15.665 21.336 /
\plot 15.665 21.336 15.627 21.323 /
\plot 15.627 21.323 15.589 21.311 /
\plot 15.589 21.311 15.551 21.298 /
\plot 15.551 21.298 15.511 21.283 /
\plot 15.511 21.283 15.469 21.266 /
\plot 15.469 21.266 15.428 21.249 /
\plot 15.428 21.249 15.386 21.232 /
\plot 15.386 21.232 15.344 21.215 /
\plot 15.344 21.215 15.301 21.196 /
\plot 15.301 21.196 15.261 21.177 /
\plot 15.261 21.177 15.221 21.156 /
\plot 15.221 21.156 15.185 21.137 /
\plot 15.185 21.137 15.149 21.118 /
\plot 15.149 21.118 15.115 21.097 /
\plot 15.115 21.097 15.081 21.076 /
\plot 15.081 21.076 15.050 21.054 /
\plot 15.050 21.054 15.016 21.033 /
\plot 15.016 21.033 14.984 21.010 /
\plot 14.984 21.010 14.950 20.985 /
\plot 14.950 20.985 14.918 20.957 /
\plot 14.918 20.957 14.884 20.930 /
\plot 14.884 20.930 14.851 20.902 /
\plot 14.851 20.902 14.817 20.872 /
\plot 14.817 20.872 14.783 20.843 /
\plot 14.783 20.843 14.751 20.813 /
\plot 14.751 20.813 14.719 20.784 /
\plot 14.719 20.784 14.690 20.754 /
\plot 14.690 20.754 14.660 20.726 /
\plot 14.660 20.726 14.633 20.699 /
\plot 14.633 20.699 14.605 20.673 /
\plot 14.605 20.673 14.580 20.648 /
\plot 14.580 20.648 14.552 20.623 /
\plot 14.552 20.623 14.525 20.597 /
\plot 14.525 20.597 14.497 20.570 /
\plot 14.497 20.570 14.470 20.544 /
\plot 14.470 20.544 14.440 20.517 /
\plot 14.440 20.517 14.410 20.487 /
\plot 14.410 20.487 14.381 20.460 /
\plot 14.381 20.460 14.351 20.430 /
\plot 14.351 20.430 14.319 20.400 /
\plot 14.319 20.400 14.290 20.371 /
\plot 14.290 20.371 14.260 20.341 /
\plot 14.260 20.341 14.230 20.314 /
\plot 14.230 20.314 14.201 20.286 /
\plot 14.201 20.286 14.173 20.259 /
\plot 14.173 20.259 14.146 20.231 /
\plot 14.146 20.231 14.118 20.204 /
\plot 14.118 20.204 14.093 20.178 /
\plot 14.093 20.178 14.069 20.155 /
\plot 14.069 20.155 14.044 20.127 /
\plot 14.044 20.127 14.017 20.100 /
\plot 14.017 20.100 13.989 20.072 /
\plot 13.989 20.072 13.962 20.043 /
\plot 13.962 20.043 13.934 20.013 /
\plot 13.934 20.013 13.904 19.981 /
\plot 13.904 19.981 13.875 19.950 /
\plot 13.875 19.950 13.847 19.918 /
\plot 13.847 19.918 13.818 19.886 /
\plot 13.818 19.886 13.790 19.854 /
\plot 13.790 19.854 13.763 19.823 /
\plot 13.763 19.823 13.735 19.793 /
\plot 13.735 19.793 13.710 19.763 /
\plot 13.710 19.763 13.684 19.734 /
\plot 13.684 19.734 13.661 19.704 /
\plot 13.661 19.704 13.638 19.674 /
\plot 13.638 19.674 13.612 19.645 /
\plot 13.612 19.645 13.589 19.615 /
\plot 13.589 19.615 13.566 19.586 /
\plot 13.566 19.586 13.542 19.554 /
\plot 13.542 19.554 13.517 19.520 /
\plot 13.517 19.520 13.492 19.486 /
\plot 13.492 19.486 13.466 19.450 /
\plot 13.466 19.450 13.441 19.416 /
\plot 13.441 19.416 13.413 19.378 /
\plot 13.413 19.378 13.388 19.342 /
\plot 13.388 19.342 13.363 19.306 /
\plot 13.363 19.306 13.337 19.270 /
\plot 13.337 19.270 13.312 19.234 /
\plot 13.312 19.234 13.286 19.198 /
\plot 13.286 19.198 13.263 19.162 /
\plot 13.263 19.162 13.240 19.128 /
\plot 13.240 19.128 13.216 19.094 /
\plot 13.216 19.094 13.193 19.061 /
\plot 13.193 19.061 13.170 19.029 /
\plot 13.170 19.029 13.149 18.999 /
\plot 13.149 18.999 13.128 18.965 /
\plot 13.128 18.965 13.104 18.934 /
\plot 13.104 18.934 13.081 18.898 /
\plot 13.081 18.898 13.056 18.862 /
\plot 13.056 18.862 13.030 18.826 /
\plot 13.030 18.826 13.005 18.788 /
\plot 13.005 18.788 12.979 18.749 /
\plot 12.979 18.749 12.952 18.709 /
\plot 12.952 18.709 12.926 18.671 /
\plot 12.926 18.671 12.899 18.631 /
\plot 12.899 18.631 12.874 18.591 /
\plot 12.874 18.591 12.848 18.553 /
\plot 12.848 18.553 12.823 18.512 /
\plot 12.823 18.512 12.797 18.474 /
\plot 12.797 18.474 12.774 18.436 /
\plot 12.774 18.436 12.751 18.398 /
\plot 12.751 18.398 12.728 18.362 /
\plot 12.728 18.362 12.706 18.324 /
\plot 12.706 18.324 12.683 18.288 /
\plot 12.683 18.288 12.660 18.250 /
\plot 12.660 18.250 12.636 18.210 /
\plot 12.636 18.210 12.613 18.169 /
\plot 12.613 18.169 12.590 18.129 /
\plot 12.590 18.129 12.565 18.087 /
\plot 12.565 18.087 12.541 18.042 /
\plot 12.541 18.042 12.516 17.998 /
\plot 12.516 17.998 12.490 17.951 /
\plot 12.490 17.951 12.467 17.905 /
\plot 12.467 17.905 12.442 17.860 /
\plot 12.442 17.860 12.418 17.814 /
\plot 12.418 17.814 12.393 17.769 /
\plot 12.393 17.769 12.372 17.723 /
\plot 12.372 17.723 12.349 17.681 /
\plot 12.349 17.681 12.327 17.638 /
\plot 12.327 17.638 12.306 17.596 /
\plot 12.306 17.596 12.287 17.556 /
\plot 12.287 17.556 12.268 17.518 /
\plot 12.268 17.518 12.251 17.477 /
\plot 12.251 17.477 12.230 17.435 /
\plot 12.230 17.435 12.211 17.393 /
\plot 12.211 17.393 12.190 17.348 /
\plot 12.190 17.348 12.171 17.304 /
\plot 12.171 17.304 12.152 17.259 /
\plot 12.152 17.259 12.131 17.213 /
\plot 12.131 17.213 12.112 17.166 /
\plot 12.112 17.166 12.093 17.120 /
\plot 12.093 17.120 12.071 17.073 /
\plot 12.071 17.073 12.054 17.029 /
\plot 12.054 17.029 12.035 16.982 /
\plot 12.035 16.982 12.018 16.940 /
\plot 12.018 16.940 12.002 16.897 /
\plot 12.002 16.897 11.987 16.857 /
\plot 11.987 16.857 11.972 16.819 /
\plot 11.972 16.819 11.957 16.781 /
\plot 11.957 16.781 11.944 16.745 /
\plot 11.944 16.745 11.934 16.711 /
\plot 11.934 16.711 11.919 16.673 /
\plot 11.919 16.673 11.906 16.635 /
\plot 11.906 16.635 11.894 16.599 /
\plot 11.894 16.599 11.881 16.561 /
\plot 11.881 16.561 11.868 16.523 /
\plot 11.868 16.523 11.858 16.482 /
\plot 11.858 16.482 11.845 16.444 /
\plot 11.845 16.444 11.834 16.406 /
\plot 11.834 16.406 11.826 16.370 /
\plot 11.826 16.370 11.815 16.334 /
\plot 11.815 16.334 11.807 16.298 /
\plot 11.807 16.298 11.800 16.264 /
\plot 11.800 16.264 11.794 16.233 /
\plot 11.794 16.233 11.788 16.203 /
\plot 11.788 16.203 11.783 16.173 /
\plot 11.783 16.173 11.779 16.144 /
\plot 11.779 16.144 11.775 16.112 /
\plot 11.775 16.112 11.773 16.080 /
\plot 11.773 16.080 11.771 16.046 /
\plot 11.771 16.046 11.769 16.013 /
\plot 11.769 16.013 11.769 15.979 /
\plot 11.769 15.979 11.771 15.943 /
\plot 11.771 15.943 11.773 15.907 /
\plot 11.773 15.907 11.777 15.869 /
\plot 11.777 15.869 11.781 15.833 /
\plot 11.781 15.833 11.786 15.799 /
\plot 11.786 15.799 11.792 15.763 /
\plot 11.792 15.763 11.800 15.729 /
\plot 11.800 15.729 11.807 15.697 /
\plot 11.807 15.697 11.817 15.663 /
\plot 11.817 15.663 11.826 15.634 /
\plot 11.826 15.634 11.834 15.604 /
\plot 11.834 15.604 11.845 15.570 /
\plot 11.845 15.570 11.858 15.536 /
\plot 11.858 15.536 11.872 15.498 /
\plot 11.872 15.498 11.889 15.456 /
\plot 11.889 15.456 11.908 15.409 /
\plot 11.908 15.409 11.930 15.359 /
\plot 11.930 15.359 11.951 15.306 /
\plot 11.951 15.306 11.972 15.255 /
\plot 11.972 15.255 11.993 15.206 /
\plot 11.993 15.206 12.010 15.166 /
\plot 12.010 15.166 12.023 15.138 /
\plot 12.023 15.138 12.029 15.121 /
\plot 12.029 15.121 12.033 15.115 /
\plot 12.033 15.115 12.033 15.113 /
}%
%
% Fig POLYLINE object
%
\linethickness= 0.500pt
\setplotsymbol ({\thinlinefont .})
\setsolid
{\putrule from 16.891 21.273 to 16.891 21.270
\plot 16.891 21.270 16.895 21.264 /
\plot 16.895 21.264 16.902 21.247 /
\plot 16.902 21.247 16.912 21.220 /
\plot 16.912 21.220 16.925 21.184 /
\plot 16.925 21.184 16.940 21.143 /
\plot 16.940 21.143 16.957 21.101 /
\plot 16.957 21.101 16.969 21.059 /
\plot 16.969 21.059 16.982 21.018 /
\plot 16.982 21.018 16.995 20.983 /
\plot 16.995 20.983 17.003 20.949 /
\plot 17.003 20.949 17.012 20.917 /
\plot 17.012 20.917 17.018 20.887 /
\plot 17.018 20.887 17.024 20.853 /
\plot 17.024 20.853 17.029 20.826 /
\plot 17.029 20.826 17.031 20.796 /
\plot 17.031 20.796 17.035 20.764 /
\plot 17.035 20.764 17.037 20.731 /
\plot 17.037 20.731 17.039 20.697 /
\plot 17.039 20.697 17.041 20.659 /
\putrule from 17.041 20.659 to 17.041 20.621
\putrule from 17.041 20.621 to 17.041 20.580
\putrule from 17.041 20.580 to 17.041 20.540
\plot 17.041 20.540 17.039 20.500 /
\plot 17.039 20.500 17.037 20.460 /
\plot 17.037 20.460 17.035 20.419 /
\plot 17.035 20.419 17.031 20.379 /
\plot 17.031 20.379 17.029 20.339 /
\plot 17.029 20.339 17.022 20.301 /
\plot 17.022 20.301 17.018 20.261 /
\plot 17.018 20.261 17.014 20.227 /
\plot 17.014 20.227 17.007 20.191 /
\plot 17.007 20.191 17.001 20.153 /
\plot 17.001 20.153 16.995 20.115 /
\plot 16.995 20.115 16.986 20.074 /
\plot 16.986 20.074 16.978 20.034 /
\plot 16.978 20.034 16.967 19.990 /
\plot 16.967 19.990 16.959 19.947 /
\plot 16.959 19.947 16.948 19.903 /
\plot 16.948 19.903 16.938 19.856 /
\plot 16.938 19.856 16.925 19.812 /
\plot 16.925 19.812 16.914 19.768 /
\plot 16.914 19.768 16.904 19.723 /
\plot 16.904 19.723 16.891 19.681 /
\plot 16.891 19.681 16.878 19.638 /
\plot 16.878 19.638 16.868 19.596 /
\plot 16.868 19.596 16.855 19.556 /
\plot 16.855 19.556 16.844 19.516 /
\plot 16.844 19.516 16.832 19.475 /
\plot 16.832 19.475 16.819 19.433 /
\plot 16.819 19.433 16.804 19.391 /
\plot 16.804 19.391 16.792 19.348 /
\plot 16.792 19.348 16.777 19.304 /
\plot 16.777 19.304 16.762 19.257 /
\plot 16.762 19.257 16.745 19.211 /
\plot 16.745 19.211 16.730 19.164 /
\plot 16.730 19.164 16.713 19.116 /
\plot 16.713 19.116 16.696 19.069 /
\plot 16.696 19.069 16.679 19.022 /
\plot 16.679 19.022 16.662 18.976 /
\plot 16.662 18.976 16.645 18.934 /
\plot 16.645 18.934 16.629 18.889 /
\plot 16.629 18.889 16.612 18.849 /
\plot 16.612 18.849 16.595 18.809 /
\plot 16.595 18.809 16.580 18.771 /
\plot 16.580 18.771 16.563 18.733 /
\plot 16.563 18.733 16.546 18.697 /
\plot 16.546 18.697 16.529 18.658 /
\plot 16.529 18.658 16.512 18.622 /
\plot 16.512 18.622 16.495 18.584 /
\plot 16.495 18.584 16.476 18.546 /
\plot 16.476 18.546 16.457 18.508 /
\plot 16.457 18.508 16.438 18.470 /
\plot 16.438 18.470 16.419 18.430 /
\plot 16.419 18.430 16.398 18.394 /
\plot 16.398 18.394 16.379 18.356 /
\plot 16.379 18.356 16.360 18.320 /
\plot 16.360 18.320 16.341 18.284 /
\plot 16.341 18.284 16.322 18.250 /
\plot 16.322 18.250 16.303 18.218 /
\plot 16.303 18.218 16.286 18.186 /
\plot 16.286 18.186 16.269 18.159 /
\plot 16.269 18.159 16.252 18.129 /
\plot 16.252 18.129 16.235 18.102 /
\plot 16.235 18.102 16.214 18.068 /
\plot 16.214 18.068 16.192 18.036 /
\plot 16.192 18.036 16.171 18.002 /
\plot 16.171 18.002 16.150 17.968 /
\plot 16.150 17.968 16.129 17.935 /
\plot 16.129 17.935 16.106 17.901 /
\plot 16.106 17.901 16.082 17.867 /
\plot 16.082 17.867 16.061 17.833 /
\plot 16.061 17.833 16.038 17.801 /
\plot 16.038 17.801 16.017 17.772 /
\plot 16.017 17.772 15.996 17.742 /
\plot 15.996 17.742 15.977 17.714 /
\plot 15.977 17.714 15.958 17.689 /
\plot 15.958 17.689 15.938 17.664 /
\plot 15.938 17.664 15.919 17.638 /
\plot 15.919 17.638 15.900 17.613 /
\plot 15.900 17.613 15.881 17.587 /
\plot 15.881 17.587 15.860 17.562 /
\plot 15.860 17.562 15.839 17.534 /
\plot 15.839 17.534 15.818 17.507 /
\plot 15.818 17.507 15.797 17.479 /
\plot 15.797 17.479 15.776 17.452 /
\plot 15.776 17.452 15.754 17.424 /
\plot 15.754 17.424 15.733 17.397 /
\plot 15.733 17.397 15.714 17.371 /
\plot 15.714 17.371 15.695 17.346 /
\plot 15.695 17.346 15.676 17.323 /
\plot 15.676 17.323 15.659 17.297 /
\plot 15.659 17.297 15.640 17.274 /
\plot 15.640 17.274 15.623 17.249 /
\plot 15.623 17.249 15.604 17.223 /
\plot 15.604 17.223 15.585 17.198 /
\plot 15.585 17.198 15.566 17.170 /
\plot 15.566 17.170 15.547 17.143 /
\plot 15.547 17.143 15.526 17.113 /
\plot 15.526 17.113 15.505 17.086 /
\plot 15.505 17.086 15.486 17.058 /
\plot 15.486 17.058 15.466 17.031 /
\plot 15.466 17.031 15.447 17.005 /
\plot 15.447 17.005 15.428 16.980 /
\plot 15.428 16.980 15.411 16.957 /
\plot 15.411 16.957 15.395 16.933 /
\plot 15.395 16.933 15.375 16.910 /
\plot 15.375 16.910 15.356 16.887 /
\plot 15.356 16.887 15.337 16.863 /
\plot 15.337 16.863 15.318 16.840 /
\plot 15.318 16.840 15.297 16.815 /
\plot 15.297 16.815 15.276 16.789 /
\plot 15.276 16.789 15.253 16.764 /
\plot 15.253 16.764 15.229 16.739 /
\plot 15.229 16.739 15.208 16.715 /
\plot 15.208 16.715 15.185 16.690 /
\plot 15.185 16.690 15.162 16.667 /
\plot 15.162 16.667 15.141 16.643 /
\plot 15.141 16.643 15.119 16.622 /
\plot 15.119 16.622 15.098 16.599 /
\plot 15.098 16.599 15.075 16.578 /
\plot 15.075 16.578 15.052 16.554 /
\plot 15.052 16.554 15.028 16.531 /
\plot 15.028 16.531 15.003 16.508 /
\plot 15.003 16.508 14.975 16.482 /
\plot 14.975 16.482 14.946 16.455 /
\plot 14.946 16.455 14.918 16.427 /
\plot 14.918 16.427 14.887 16.402 /
\plot 14.887 16.402 14.857 16.375 /
\plot 14.857 16.375 14.827 16.347 /
\plot 14.827 16.347 14.798 16.322 /
\plot 14.798 16.322 14.768 16.296 /
\plot 14.768 16.296 14.740 16.271 /
\plot 14.740 16.271 14.711 16.245 /
\plot 14.711 16.245 14.685 16.224 /
\plot 14.685 16.224 14.658 16.201 /
\plot 14.658 16.201 14.630 16.178 /
\plot 14.630 16.178 14.601 16.152 /
\plot 14.601 16.152 14.571 16.129 /
\plot 14.571 16.129 14.539 16.101 /
\plot 14.539 16.101 14.508 16.076 /
\plot 14.508 16.076 14.474 16.049 /
\plot 14.474 16.049 14.440 16.021 /
\plot 14.440 16.021 14.408 15.996 /
\plot 14.408 15.996 14.374 15.968 /
\plot 14.374 15.968 14.343 15.943 /
\plot 14.343 15.943 14.311 15.919 /
\plot 14.311 15.919 14.279 15.894 /
\plot 14.279 15.894 14.249 15.871 /
\plot 14.249 15.871 14.220 15.847 /
\plot 14.220 15.847 14.188 15.826 /
\plot 14.188 15.826 14.158 15.803 /
\plot 14.188 15.826 14.158 15.803 /
\plot 14.158 15.803 14.127 15.780 /
\plot 14.127 15.780 14.093 15.754 /
\plot 14.093 15.754 14.059 15.731 /
\plot 14.059 15.731 14.023 15.706 /
\plot 14.023 15.706 13.987 15.680 /
\plot 13.987 15.680 13.951 15.655 /
\plot 13.951 15.655 13.915 15.629 /
\plot 13.915 15.629 13.877 15.604 /
\plot 13.877 15.604 13.841 15.579 /
\plot 13.841 15.579 13.805 15.555 /
\plot 13.805 15.555 13.771 15.532 /
\plot 13.771 15.532 13.735 15.511 /
\plot 13.735 15.511 13.701 15.488 /
\plot 13.701 15.488 13.669 15.466 /
\plot 13.669 15.466 13.633 15.447 /
\plot 13.633 15.447 13.600 15.426 /
\plot 13.600 15.426 13.564 15.405 /
\plot 13.564 15.405 13.528 15.382 /
\plot 13.528 15.382 13.490 15.361 /
\plot 13.490 15.361 13.451 15.339 /
\plot 13.451 15.339 13.411 15.316 /
\plot 13.411 15.316 13.371 15.295 /
\plot 13.371 15.295 13.333 15.274 /
\plot 13.333 15.274 13.295 15.255 /
\plot 13.295 15.255 13.259 15.236 /
\plot 13.259 15.236 13.223 15.217 /
\plot 13.223 15.217 13.191 15.202 /
\plot 13.191 15.202 13.159 15.187 /
\plot 13.159 15.187 13.130 15.172 /
\plot 13.130 15.172 13.102 15.160 /
\plot 13.102 15.160 13.066 15.145 /
\plot 13.066 15.145 13.032 15.132 /
\plot 13.032 15.132 12.998 15.119 /
\plot 12.998 15.119 12.965 15.109 /
\plot 12.965 15.109 12.931 15.098 /
\plot 12.931 15.098 12.897 15.088 /
\plot 12.897 15.088 12.865 15.079 /
\plot 12.865 15.079 12.835 15.073 /
\plot 12.835 15.073 12.806 15.066 /
\plot 12.806 15.066 12.778 15.060 /
\plot 12.778 15.060 12.753 15.054 /
\plot 12.753 15.054 12.728 15.050 /
\plot 12.728 15.050 12.700 15.045 /
\plot 12.700 15.045 12.675 15.039 /
\plot 12.675 15.039 12.647 15.035 /
\plot 12.647 15.035 12.617 15.030 /
\plot 12.617 15.030 12.590 15.026 /
\plot 12.590 15.026 12.560 15.022 /
\plot 12.560 15.022 12.531 15.020 /
\plot 12.531 15.020 12.501 15.018 /
\plot 12.501 15.018 12.476 15.016 /
\putrule from 12.476 15.016 to 12.450 15.016
\putrule from 12.450 15.016 to 12.427 15.016
\plot 12.427 15.016 12.404 15.018 /
\plot 12.404 15.018 12.378 15.022 /
\plot 12.378 15.022 12.353 15.028 /
\plot 12.353 15.028 12.325 15.037 /
\plot 12.325 15.037 12.294 15.047 /
\plot 12.294 15.047 12.260 15.062 /
\plot 12.260 15.062 12.226 15.079 /
\plot 12.226 15.079 12.194 15.096 /
\plot 12.194 15.096 12.173 15.107 /
\plot 12.173 15.107 12.162 15.113 /
\putrule from 12.162 15.113 to 12.160 15.113
}%
%
% Fig POLYLINE object
%
\linethickness= 0.500pt
\setplotsymbol ({\thinlinefont .})
{\plot 16.897 21.294 16.891 21.296 /
\plot 16.891 21.296 16.878 21.302 /
\plot 16.878 21.302 16.855 21.311 /
\plot 16.855 21.311 16.821 21.325 /
\plot 16.821 21.325 16.775 21.344 /
\plot 16.775 21.344 16.720 21.368 /
\plot 16.720 21.368 16.658 21.393 /
\plot 16.658 21.393 16.590 21.421 /
\plot 16.590 21.421 16.521 21.448 /
\plot 16.521 21.448 16.453 21.476 /
\plot 16.453 21.476 16.387 21.503 /
\plot 16.387 21.503 16.324 21.529 /
\plot 16.324 21.529 16.264 21.552 /
\plot 16.264 21.552 16.209 21.573 /
\plot 16.209 21.573 16.159 21.592 /
\plot 16.159 21.592 16.112 21.609 /
\plot 16.112 21.609 16.070 21.624 /
\plot 16.070 21.624 16.030 21.639 /
\plot 16.030 21.639 15.991 21.651 /
\plot 15.991 21.651 15.953 21.664 /
\plot 15.953 21.664 15.919 21.675 /
\plot 15.919 21.675 15.875 21.687 /
\plot 15.875 21.687 15.835 21.698 /
\plot 15.835 21.698 15.792 21.709 /
\plot 15.792 21.709 15.750 21.719 /
\plot 15.750 21.719 15.708 21.730 /
\plot 15.708 21.730 15.665 21.738 /
\plot 15.665 21.738 15.623 21.745 /
\plot 15.623 21.745 15.581 21.753 /
\plot 15.581 21.753 15.538 21.757 /
\plot 15.538 21.757 15.498 21.764 /
\plot 15.498 21.764 15.458 21.768 /
\plot 15.458 21.768 15.418 21.772 /
\plot 15.418 21.772 15.378 21.774 /
\plot 15.378 21.774 15.342 21.776 /
\plot 15.342 21.776 15.303 21.778 /
\plot 15.303 21.778 15.268 21.780 /
\putrule from 15.268 21.780 to 15.234 21.780
\putrule from 15.234 21.780 to 15.200 21.780
\putrule from 15.200 21.780 to 15.164 21.780
\putrule from 15.164 21.780 to 15.128 21.780
\plot 15.128 21.780 15.092 21.778 /
\plot 15.092 21.778 15.054 21.776 /
\putrule from 15.054 21.776 to 15.014 21.776
\plot 15.014 21.776 14.973 21.772 /
\plot 14.973 21.772 14.931 21.770 /
\plot 14.931 21.770 14.887 21.766 /
\plot 14.887 21.766 14.842 21.761 /
\plot 14.842 21.761 14.798 21.755 /
\plot 14.798 21.755 14.753 21.751 /
\plot 14.753 21.751 14.707 21.745 /
\plot 14.707 21.745 14.662 21.738 /
\plot 14.662 21.738 14.620 21.730 /
\plot 14.620 21.730 14.575 21.721 /
\plot 14.575 21.721 14.533 21.715 /
\plot 14.533 21.715 14.491 21.704 /
\plot 14.491 21.704 14.448 21.696 /
\plot 14.448 21.696 14.408 21.687 /
\plot 14.408 21.687 14.368 21.677 /
\plot 14.368 21.677 14.328 21.666 /
\plot 14.328 21.666 14.285 21.656 /
\plot 14.285 21.656 14.241 21.643 /
\plot 14.241 21.643 14.196 21.630 /
\plot 14.196 21.630 14.150 21.615 /
\plot 14.150 21.615 14.103 21.601 /
\plot 14.103 21.601 14.055 21.586 /
\plot 14.055 21.586 14.008 21.569 /
\plot 14.008 21.569 13.959 21.552 /
\plot 13.959 21.552 13.911 21.535 /
\plot 13.911 21.535 13.864 21.518 /
\plot 13.864 21.518 13.818 21.499 /
\plot 13.818 21.499 13.771 21.482 /
\plot 13.771 21.482 13.729 21.463 /
\plot 13.729 21.463 13.684 21.444 /
\plot 13.684 21.444 13.644 21.427 /
\plot 13.644 21.427 13.604 21.408 /
\plot 13.604 21.408 13.564 21.389 /
\plot 13.564 21.389 13.526 21.370 /
\plot 13.526 21.370 13.485 21.351 /
\plot 13.485 21.351 13.445 21.330 /
\plot 13.445 21.330 13.407 21.308 /
\plot 13.407 21.308 13.365 21.285 /
\plot 13.365 21.285 13.324 21.262 /
\plot 13.324 21.262 13.284 21.239 /
\plot 13.284 21.239 13.242 21.213 /
\plot 13.242 21.213 13.200 21.188 /
\plot 13.200 21.188 13.159 21.162 /
\plot 13.159 21.162 13.117 21.137 /
\plot 13.117 21.137 13.077 21.110 /
\plot 13.077 21.110 13.039 21.084 /
\plot 13.039 21.084 13.001 21.059 /
\plot 13.001 21.059 12.962 21.033 /
\plot 12.962 21.033 12.929 21.008 /
\plot 12.929 21.008 12.893 20.983 /
\plot 12.893 20.983 12.861 20.959 /
\plot 12.861 20.959 12.827 20.936 /
\plot 12.827 20.936 12.797 20.913 /
\plot 12.797 20.913 12.761 20.887 /
\plot 12.761 20.887 12.728 20.860 /
\plot 12.728 20.860 12.692 20.834 /
\plot 12.692 20.834 12.656 20.805 /
\plot 12.656 20.805 12.620 20.777 /
\plot 12.620 20.777 12.584 20.748 /
\plot 12.584 20.748 12.548 20.716 /
\plot 12.548 20.716 12.510 20.686 /
\plot 12.510 20.686 12.471 20.654 /
\plot 12.471 20.654 12.435 20.621 /
\plot 12.435 20.621 12.397 20.589 /
\plot 12.397 20.589 12.361 20.557 /
\plot 12.361 20.557 12.327 20.525 /
\plot 12.327 20.525 12.294 20.494 /
\plot 12.294 20.494 12.260 20.462 /
\plot 12.260 20.462 12.228 20.430 /
\plot 12.228 20.430 12.198 20.398 /
\plot 12.198 20.398 12.167 20.367 /
\plot 12.167 20.367 12.139 20.339 /
\plot 12.139 20.339 12.112 20.309 /
\plot 12.112 20.309 12.082 20.280 /
\plot 12.082 20.280 12.052 20.248 /
\plot 12.052 20.248 12.023 20.214 /
\plot 12.023 20.214 11.993 20.180 /
\plot 11.993 20.180 11.961 20.144 /
\plot 11.961 20.144 11.932 20.106 /
\plot 11.932 20.106 11.900 20.068 /
\plot 11.900 20.068 11.868 20.030 /
\plot 11.868 20.030 11.836 19.990 /
\plot 11.836 19.990 11.805 19.950 /
\plot 11.805 19.950 11.773 19.909 /
\plot 11.773 19.909 11.743 19.869 /
\plot 11.743 19.869 11.714 19.831 /
\plot 11.714 19.831 11.686 19.791 /
\plot 11.686 19.791 11.659 19.753 /
\plot 11.659 19.753 11.631 19.715 /
\plot 11.631 19.715 11.606 19.677 /
\plot 11.606 19.677 11.580 19.636 /
\plot 11.580 19.636 11.557 19.602 /
\plot 11.557 19.602 11.534 19.566 /
\plot 11.534 19.566 11.510 19.528 /
\plot 11.510 19.528 11.487 19.490 /
\plot 11.487 19.490 11.462 19.452 /
\plot 11.462 19.452 11.438 19.410 /
\plot 11.438 19.410 11.413 19.367 /
\plot 11.413 19.367 11.390 19.325 /
\plot 11.390 19.325 11.364 19.281 /
\plot 11.364 19.281 11.341 19.236 /
\plot 11.341 19.236 11.316 19.190 /
\plot 11.316 19.190 11.292 19.143 /
\plot 11.292 19.143 11.269 19.097 /
\plot 11.269 19.097 11.246 19.050 /
\plot 11.246 19.050 11.223 19.003 /
\plot 11.223 19.003 11.201 18.959 /
\plot 11.201 18.959 11.180 18.912 /
\plot 11.180 18.912 11.161 18.868 /
\plot 11.161 18.868 11.140 18.826 /
\plot 11.140 18.826 11.123 18.781 /
\plot 11.123 18.781 11.104 18.739 /
\plot 11.104 18.739 11.087 18.694 /
\plot 11.087 18.694 11.070 18.652 /
\plot 11.070 18.652 11.053 18.608 /
\plot 11.053 18.608 11.036 18.563 /
\plot 11.036 18.563 11.019 18.519 /
\plot 11.019 18.519 11.002 18.472 /
\plot 11.002 18.472 10.988 18.423 /
\plot 10.988 18.423 10.971 18.375 /
\plot 10.971 18.375 10.954 18.326 /
\plot 10.954 18.326 10.939 18.275 /
\plot 10.939 18.275 10.922 18.227 /
\plot 10.922 18.227 10.907 18.176 /
\plot 10.907 18.176 10.892 18.127 /
\plot 10.892 18.127 10.880 18.078 /
\plot 10.880 18.078 10.865 18.032 /
\plot 10.865 18.032 10.852 17.985 /
\plot 10.852 17.985 10.842 17.941 /
\plot 10.842 17.941 10.829 17.899 /
\plot 10.829 17.899 10.820 17.858 /
\plot 10.820 17.858 10.810 17.818 /
\plot 10.810 17.818 10.801 17.782 /
\plot 10.801 17.782 10.793 17.746 /
\plot 10.793 17.746 10.787 17.710 /
\plot 10.787 17.710 10.776 17.666 /
\plot 10.776 17.666 10.767 17.621 /
\plot 10.767 17.621 10.759 17.579 /
\plot 10.759 17.579 10.753 17.534 /
\plot 10.753 17.534 10.746 17.492 /
\plot 10.746 17.492 10.742 17.448 /
\plot 10.742 17.448 10.738 17.405 /
\plot 10.738 17.405 10.736 17.363 /
\plot 10.736 17.363 10.734 17.321 /
\putrule from 10.734 17.321 to 10.734 17.280
\putrule from 10.734 17.280 to 10.734 17.240
\plot 10.734 17.240 10.736 17.200 /
\plot 10.736 17.200 10.740 17.162 /
\plot 10.740 17.162 10.742 17.126 /
\plot 10.742 17.126 10.748 17.088 /
\plot 10.748 17.088 10.755 17.050 /
\plot 10.755 17.050 10.759 17.016 /
\plot 10.759 17.016 10.767 16.980 /
\plot 10.767 16.980 10.774 16.942 /
\plot 10.774 16.942 10.782 16.904 /
\plot 10.782 16.904 10.793 16.863 /
\plot 10.793 16.863 10.803 16.821 /
\plot 10.803 16.821 10.816 16.777 /
\plot 10.816 16.777 10.827 16.732 /
\plot 10.827 16.732 10.839 16.688 /
\plot 10.839 16.688 10.854 16.643 /
\plot 10.854 16.643 10.867 16.599 /
\plot 10.867 16.599 10.882 16.554 /
\plot 10.882 16.554 10.897 16.512 /
\plot 10.897 16.512 10.909 16.472 /
\plot 10.909 16.472 10.924 16.432 /
\plot 10.924 16.432 10.939 16.394 /
\plot 10.939 16.394 10.952 16.355 /
\plot 10.952 16.355 10.966 16.320 /
\plot 10.966 16.320 10.981 16.279 /
\plot 10.981 16.279 10.996 16.241 /
\plot 10.996 16.241 11.013 16.201 /
\plot 11.013 16.201 11.030 16.161 /
\plot 11.030 16.161 11.049 16.121 /
\plot 11.049 16.121 11.068 16.080 /
\plot 11.068 16.080 11.087 16.040 /
\plot 11.087 16.040 11.108 16.000 /
\plot 11.108 16.000 11.127 15.962 /
\plot 11.127 15.962 11.148 15.924 /
\plot 11.148 15.924 11.170 15.886 /
\plot 11.170 15.886 11.191 15.852 /
\plot 11.191 15.852 11.210 15.818 /
\plot 11.210 15.818 11.231 15.786 /
\plot 11.231 15.786 11.252 15.756 /
\plot 11.252 15.756 11.273 15.727 /
\plot 11.273 15.727 11.295 15.697 /
\plot 11.295 15.697 11.316 15.670 /
\plot 11.316 15.670 11.339 15.640 /
\plot 11.339 15.640 11.364 15.610 /
\plot 11.364 15.610 11.390 15.581 /
\plot 11.390 15.581 11.417 15.551 /
\plot 11.417 15.551 11.447 15.519 /
\plot 11.447 15.519 11.477 15.490 /
\plot 11.477 15.490 11.506 15.462 /
\plot 11.506 15.462 11.538 15.435 /
\plot 11.538 15.435 11.568 15.407 /
\plot 11.568 15.407 11.599 15.382 /
\plot 11.599 15.382 11.631 15.356 /
\plot 11.631 15.356 11.661 15.335 /
\plot 11.661 15.335 11.692 15.314 /
\plot 11.692 15.314 11.722 15.293 /
\plot 11.722 15.293 11.750 15.276 /
\plot 11.750 15.276 11.779 15.259 /
\plot 11.779 15.259 11.809 15.240 /
\plot 11.809 15.240 11.841 15.223 /
\plot 11.841 15.223 11.872 15.206 /
\plot 11.872 15.206 11.906 15.189 /
\plot 11.906 15.189 11.942 15.172 /
\plot 11.942 15.172 11.978 15.153 /
\plot 11.978 15.153 12.014 15.136 /
\plot 12.014 15.136 12.052 15.119 /
\plot 12.052 15.119 12.090 15.102 /
\plot 12.090 15.102 12.129 15.085 /
\plot 12.129 15.085 12.167 15.069 /
\plot 12.167 15.069 12.203 15.054 /
\plot 12.203 15.054 12.239 15.037 /
\plot 12.239 15.037 12.275 15.022 /
\plot 12.275 15.022 12.311 15.007 /
\plot 12.311 15.007 12.347 14.990 /
\plot 12.347 14.990 12.383 14.975 /
\plot 12.383 14.975 12.418 14.961 /
\plot 12.418 14.961 12.457 14.944 /
\plot 12.457 14.944 12.497 14.927 /
\plot 12.497 14.927 12.537 14.910 /
\plot 12.537 14.910 12.577 14.891 /
\plot 12.577 14.891 12.622 14.874 /
\plot 12.622 14.874 12.664 14.855 /
\plot 12.664 14.855 12.708 14.838 /
\plot 12.708 14.838 12.753 14.821 /
\plot 12.753 14.821 12.797 14.804 /
\plot 12.797 14.804 12.842 14.787 /
\plot 12.842 14.787 12.886 14.772 /
\plot 12.886 14.772 12.929 14.757 /
\plot 12.929 14.757 12.969 14.745 /
\plot 12.969 14.745 13.009 14.732 /
\plot 13.009 14.732 13.049 14.721 /
\plot 13.049 14.721 13.087 14.711 /
\plot 13.087 14.711 13.125 14.702 /
\plot 13.125 14.702 13.166 14.694 /
\plot 13.166 14.694 13.206 14.685 /
\plot 13.206 14.685 13.248 14.679 /
\plot 13.248 14.679 13.295 14.673 /
\plot 13.295 14.673 13.343 14.666 /
\plot 13.343 14.666 13.399 14.660 /
\plot 13.399 14.660 13.456 14.654 /
\plot 13.456 14.654 13.517 14.649 /
\plot 13.517 14.649 13.581 14.643 /
\plot 13.581 14.643 13.644 14.639 /
\plot 13.644 14.639 13.703 14.635 /
\plot 13.703 14.635 13.754 14.633 /
\plot 13.754 14.633 13.796 14.628 /
\putrule from 13.796 14.628 to 13.824 14.628
\plot 13.824 14.628 13.841 14.626 /
\putrule from 13.841 14.626 to 13.849 14.626
}%
%
% Fig TEXT object
%
\put{$\nu$} [lB] at 15.303 22.447
%
% Fig TEXT object
%
\put{$S$} [lB] at 12.986 21.495
%
% Fig TEXT object
%
%
\put{$e_{n+1}$} [lB] at 11.190 17.479
%
% Fig TEXT object
%
\put{$0$} [lB] at 12.033 14.319
%
% Fig TEXT object
%
\put{$p$} [lB] at 13.684 13.926
%
% Fig TEXT object
%
\put{$\Gamma$} [lB] at 16.415 17.653
%
% Fig TEXT object
%
\put{$e_1$} [lB] at 15.367 14.986
%
% Fig TEXT object
%
\put{$^{_\bullet}$} [lB] at 13.800 14.446
\linethickness=0pt
\putrectangle corners at 10.708 22.638 and 17.067 13.906
\endpicture}
$$

\bigskip
Let $\Cal E_0=R^{-\frac 1{2n}} \Cal L\cap\Gamma_0$ with $\Gamma_0$ the intersection of $\Gamma$ with a fixed neighborhood of $0$ and $\Cal L$ a suitable regular lattice
in $V$ of volume $O(1)$.

Proceeding as in \S4, the lattice $\Cal L$ can be chosen s.t.
$$
\pi_{x'}(\Cal L^*)+[-CR^{\frac {n-1}{2n}}, CR ^{\frac {n-1}{2n}}]e_1\tag 5.7
$$
is $O(R^{-\frac 1{2n}})$ dense in $[e_1, \ldots, e_n]$, hence
$$
\pi_{x'}(R^{\frac 1{2n}} \Cal L^*+[-CR^{\frac 12}, CR^{\frac 12}]\nu)\tag 5.8
$$
is $o(1)$-dense on $[e_1, \ldots, e_n]$.

Note that for $n=2$, a suitable choice of the curve $\Gamma$ permits to ensure that
$$
|\Cal E_0| =|\Cal L\cap R^{\frac 14} \Gamma_0|\sim R^{\frac 18}.\tag 5.9
$$
For $n\geq 3$, taking $\Gamma$ to be an $(n-1)$-sphere, one may for appropriate $R$ obtain
$$
|\Cal E_0| =|\Cal L\cap R^{\frac 1{2n}} \Gamma_0|\sim R^{\frac {n-2}{2n}}.\tag 5.10
$$
Let $\rho>0$ be a small constant.
Recalling that $\nu$ is tangent to $S$ at the points of $\Gamma$, each point from the set $D = \Cal E_0+\big[0, \frac \rho 
{\sqrt R}\big]\nu$ is at distance at most
$0\big(\frac {\rho^2}{R}\big)$ from $S$.
This clearly allows to construct a measure $\mu$ on $S$, $\frac {d\mu}{d\sigma} \in L^2(S, \sigma), \big\Vert\frac {d\mu}{d\sigma}\Vert_2=1$, such that for
$|x|<R$,
$$
\hat\mu(x) \sim\Big(\frac \rho{\sqrt R}\Big)^{-\frac 12} \Big(\frac\rho{R}\Big)^{\frac {n-1}2} |\Cal E_0|^{-\frac 12} \sum_{\xi\in\Cal E_0}
\Big[e(x, \xi) \Big(\int_0^{\frac \rho{\sqrt R}}
 e(sx.\nu)ds\Big)+O\Big(\frac {\rho^2}{\sqrt R}\Big)\Big].\tag 5.11
$$
If $x=y+\lambda\nu$ with $y\in R^{\frac 1{2n}}\Cal L^*\subset V$ and $|\lambda|< CR^{\frac 12}$, we have $\Vert x.\xi\Vert =\Vert y.\xi\Vert=0$
($\Vert \ \Vert$ referring to the distance to the nearest integer) by definition of $\Cal E_0$, while
$\Vert sx.\nu\Vert\leq \frac \rho{\sqrt R} |\lambda|< O(\rho)$.
Hence
$$
(5.11) \sim R^{-\frac {n-1}2 -\frac 14} |\Cal E_0|^{\frac 12}\tag 5.12
$$
holds for $x$ in an $o(1)$-neighborhood of $R^{\frac 1{2n}}\Cal L^* +[-CR^{\frac 12}, CR^{\frac 12}]\nu$.

Recalling (5.8), the preceding implies that $\mu$ introduced above satisfies
$$
\sup_{|x_{n+1}|\lesssim R} |\hat\mu(x', x_{n+1})|\gtrsim R^{-\frac n2 +\frac 14}|\Cal E_0|^{\frac 12}\tag 5.13
$$
for $x' \in B_R\cap \Bbb R^n$, and therefore
$$
\Vert(\hat\mu)^*\Vert_{L^2_{[|x'|\leq R]}} \gtrsim R^{\frac 14} |\Cal E_0|^{\frac 12}.\tag 5.14
$$
For $n=2$, (5.9) gives a lower bound $R^{\frac 5{16}}$ while for $n\geq 3$, we obtain
$R^{\frac 12 -\frac 1{2n}}$.


\Refs
\widestnumber\no{XXXXXX}

\ref
\no{[B1]}
\by J.~Bourgain
\paper A remark on Schr\"odinger operators
\jour Israel J.~Math 77 (1992), 1--16
\endref

\ref 
\no{[B2]}
\by J.~Bourgain
\paper
Some new estimates on oscillatory integrals
\jour Essays on Fourier Analysis in Honor of Elias M.~Stein (Princeton, NJ 1991), Princeton Math.Ser., Vol. 42, Princeton University Press,
New Jersey, 1995, pp. 83--112
\endref

\ref
\no{[B-G]}
\by J.~Bourgain, L.~Guth
\paper
Bounds on oscillatory integral operators based on multilinear estimates
\jour GAFA, Vol. 21, No 6, 1239--1295 (2011)
\endref

\ref
\no{[B-C-T]}
\by J.~Bennett, T.~Carbery, T.~Tao
\paper On the multilinear restriction and Kakeya conjectures
\jour Acta Math. 196 (2006), 261--302
\endref

\ref
\no {[C]}
\by L.~Carleson
\paper Some analytic problems related to statistical mechanics
\jour Euclidean Harmonic Analysis (Proc. Sem., Univ. Maryland, College Park, Md. 1979), Lecture Notes in Math., Vol. 779, Springer, Berlin,
1980, pp. 5--45
\endref

\ref
\no{[D-K]} \by B.~E.J.~Dahlberg, C.~E.~Kenig
\paper A note on the almost everywhere behavior of solutions to the Schr\"odinger equation
\jour Harmonic Analysis (Minneapolis, Minn, 1981), Lecture Notes in Math., Vol. 908, Springer, Berlin, 1982, pp. 205--209
\endref

\ref
\no{[L]}\by S.~Lee
\paper
On pointwise convergence of the solutions to Schr\"odinger equations in $\Bbb R^2$
\jour IMRN, Vol 2006, 1--21
\endref

\ref
\no {[M-V-V]}
\by A.~Moyua, A.~Vargas, L.~Vega
\paper Schr\"odinger maximal function and restriction properties of the Fourier transform
\jour International Mathematics Research Notices 1996, no.16, 793--815
\endref

\ref\no{[S]}
\by P.~Sj\"olin
\paper Regularity of solutions to the Schr\"odinger equation
\jour Duke Mathematical Journal 55 (1987), no. 3, 699--715
\endref

\ref\no{[T]}
\by T.~Tao
\paper A sharp bilinear restrictions estimate for paraboloids
\jour Geometric and Functional Analysis 13 (2003), no. 6, 1359--1384
\endref

\ref\no{[T-V]}
\by T.~Tao, A.~Vargas
\paper
A bilinear approach to cone multipliers. I. Restriction estimates
\jour Geometric and Functional Analysis 10 (2000), no. 1, 185--215
\endref

\ref\no{[T-V2]}
\by T.~Tao, A.~Vargas
\paper A bilinear approach to cone multipliers. II. Applications
\jour Geometric and Functional Analysis 10 (2000), no. 1, 216--258
\endref


\ref\no
{[V]}\by L.~Vega
\paper
Schr\"odinger equations: pointwise convergence to the initial data
\jour Proceedings of the American Mathematical Society 102 (1988), no. 4, 874--878
\endref
\endRefs
\enddocument
\end